\documentclass[11pt]{article}

\def\mathbb{\Bbb}
\usepackage{amsfonts,latexsym,amstext}
\usepackage{amsmath}
\usepackage[active]{srcltx}
\usepackage{amssymb}

\newtheorem{theorem}{Theorem}[section]
\newtheorem{lemma}[theorem]{Lemma}
\newtheorem{proposition}[theorem]{Proposition}
\newtheorem{definition}{Definition}[section]

\newtheorem{hypothesis}[theorem]{Hypothesis}

\newtheorem{remark}[theorem]{Remark}
\newtheorem{corollary}[theorem]{Corollary}
\numberwithin{equation}{section}
\def\qed{{\hfill\hbox{\enspace${ \square}$}} \smallskip}
\def\sqr#1#2{{\vcenter{\vbox{\hrule height .#2pt \hbox{\vrule
 width .#2pt height#1pt \kern#1pt \vrule
width .#2pt} \hrule height .#2pt}}}}
\def\square{\mathchoice\sqr54\sqr54\sqr{4.1}3\sqr{3.5}3}

\def\ds{\begin{displaystyle}}
\def\eds{\end{displaystyle}}
\def\dis{\displaystyle }
\def\<{\langle }
\def\>{\rangle }

\def\dim{\noindent \hbox{{\bf Proof.} }}
\def\R{\mathbb R}
\def\N{\mathbb N}
\def\Z{\mathbb Z}

\def\E{\mathbb E}
\def\P{\mathbb P}
\def\Q{\mathbb Q}

\def\cald{{\cal D}}
\def\cale{{\cal E}}
\def\calf{{\cal F}}
\def\calg{{\cal G}}

\def\caln{{\cal N}}

\def\call{{\cal L}}

\def\calo{{\cal O}}

\begin{document}

\title{Hamilton Jacobi Bellman equations in infinite dimensions
with quadratic and superquadratic Hamiltonian}
\date{}
 \author{
 Federica Masiero\\
 Dipartimento di Matematica e Applicazioni, Universit\`a di Milano Bicocca\\
 via Cozzi 53, 20125 Milano, Italy\\
 e-mail: federica.masiero@unimib.it,
 %gianmario.tessitore@unimib.it
}

\maketitle  
\begin{abstract}
We consider Hamilton Jacobi Bellman equations in an inifinite dimensional Hilbert space,
with quadratic (respectively superquadratic) hamiltonian and with continuous (respectively lipschitz continuous)
final conditions. This allows to study stochastic optimal control problems for suitable
controlled Ornstein Uhlenbeck process with unbounded control processes.
\end{abstract}

\section{Introduction}

In this paper we study semilinear Kolmogorov equations
in an infinite dimensional Hilbert space $H$, in particular Hamilton Jacobi Bellman equations.
More precisely, let us consider the following equation
\begin{equation}
\left\{
\begin{array}
[c]{l}%
\frac{\partial v}{\partial t}(t,x)=-\mathcal{L}v\left(  t,x\right)
+\psi\left(  \nabla v\left(  t,x\right)\sqrt{Q}  \right)+l(x)  ,\text{ \ \ \ \ }t\in\left[  0,T\right]
,\text{ }x\in H\\
v(T,x)=\phi\left(  x\right),
\end{array}
\right.  \label{Kolmo intro}%
\end{equation}
where $\mathcal L $ is the generator of the transition
semigroup $P_t$ related to the following Ornstein-Uhlenbeck process
\begin{equation}
\left\{
\begin{array}
[c]{l}%
dX_t  =AX_t dt+\sqrt{Q}dW_t
,\text{ \ \ \ }t\in\left[ 0,T\right] \\
X_0 =x,
\end{array}
\right.  \label{ornstein intro}%
\end{equation}
that is, at least formally, 
$$
(\call f)(x)=\frac{1}{2}(Tr Q \nabla^2 f)(x)+\<Ax,\nabla f(x)\>.
$$
The aim of this paper is to to consider the case
where $\psi$ has quadratic or superquadratic growth, and to apply our results
to suitable stochastic optimal control problems: to this aim we
make some regularizing assumptions on the Ornstein Uhlenbeck transition semigroup.
At first in equation (\ref{Kolmo intro}) we consider the case of final condition $\phi$
lipschitz continuous: with this assumption we can solve the Kolmogorov equation 
with $\psi$ quadratic and superquadratic. In the case of quadratic hamiltonian
we can solve equation (\ref{Kolmo intro}) also in the case of final condiotion
$\phi$ only bounded and continuous.
A similar result,  with $\psi$ quadratic and superquadratic and with
final condition $\phi$ lipschitz continuous,
is proved in \cite{Go2} by means of a detailed study
on weakly continuous semigroups, and making the assumption
that the transition semigroup $P_t$ is strong Feller.
Here we include the degenerate case and we exploit
the connection between PDEs and backward
stochastic differential equations (BSDEs in the following).

Coming into more details, we assume that $A$ and $Q$ in equation (\ref{ornstein intro})
commute, so that, see \cite{Mas}, the transition semigroup $P_t$ satisfies the following regularizing
property: for every $\phi\in C_{b}\left(  H\right)  $, for
every $\xi\in H,$ the function $P_t\phi$ is G\^ateaux differentiable in the direction
$\sqrt{Q}\xi$ and for $0<t\leq T$,%
\begin{equation}
\left|  \nabla P_{t}\left[  \phi\right]
\left(  x\right)\sqrt{Q}  \xi\right|  \leq\frac{c}{t  ^{1/2}}\left\|  \phi\right\|  _{\infty}\left|  \xi\right|  .\label{ipotesi H intro}%
\end{equation}
In order to prove existence 
and uniqueness of a mild solution $v$ of equation (\ref{Kolmo intro}),
we use the fact that $v$ can be represented in terms of the solution of a suitable
forward-backward system (FBSDE in the following):
\begin{equation}\label{fbsde intro}
    \left\{\begin{array}{l}\dis dX_\tau =
AX_\tau d\tau+ \sqrt{Q} dW_\tau,\quad \tau\in
[t,T]\subset [0,T],
\\\dis
X_t=x,
\\\dis
 dY_\tau=-\psi(Z_\tau)\;d\tau-l(X_\tau)\;d\tau+Z_\tau\;dW_\tau,
  \\\dis
  Y_T=\phi(X_T),
\end{array}\right.
\end{equation}
It is well known, see e.g. \cite{PaPe} for the finite dimensional case and \cite{fute}
for the generalization to the infinite dimensional case, that $v(t,x)=Y_t^{t,x}$,
so that estimates on $v$ can be achieved by studying the BSDE
\begin{equation}\label{bsde intro}
    \left\{\begin{array}{l}\dis
 dY_\tau=-\psi(Z_\tau)\;d\tau-l(X_\tau)\;d\tau+Z_\tau\;dW_\tau,
 \qquad \tau\in [0,T],
  \\\dis
  Y_T=\phi(X_T),
\end{array}\right.
\end{equation}
Moreover, if $\psi$ is quadratic, we can remove the lipschitz continuous assumption on $\phi$
and prove existence and uniqueness of a mild solution of equation (\ref{Kolmo intro})
with $\phi$ continuous and bounded.
The fundamental tool is an apriori estimate on $Z$,
and the classical identification $Z_\tau^{t,x}=\nabla v(\tau,X_\tau^{t,x})\sqrt{Q}$:
the fact that $A$ and $Q$ commute is crucial in proving this estimates on $\nabla v(t,x)\sqrt{Q} $
by means of backward stochastic differential equations.
This estimate is obtained with techniques similar to the ones introduced in \cite{BaoDeHu}, and specialized in \cite{Ri} in the quadratic case, 
to treat BSDEs with generator $\psi$ with superquadratic growth and in a markovian framework.
In \cite{BaoDeHu} the Markov process $X$ solves a finite dimensional stochastic differential equation,
with constant diffusion coefficient and with drift not necessarily linear as in our case.
In order to obtain an estimate on $Z_\tau^{t,x}$, some non degeneracy assumptions on the coefficients are made.
In the present paper the process $X$ is infinite dimensional and we need the coefficient $A$
and $\sqrt{Q}$ commute. We note that not in \cite{BaoDeHu} nor in \cite{Ri} the estimate on $Z$
is used in order to solve a PDE related.

\noindent We also cite the paper \cite{BriFu} where infinite dimensional
Hamilton Jacobi Bellman equations with quadratic hamiltonian are solved:
the generator $\call$ is related to a more general Markov process $X$ then the one considered here in
(\ref{ornstein intro}), and no assumptions on the coefficicent are made,
but only the case of final condition $\phi$ G\^ateaux differentiable is treated.

We apply these results on equation (\ref{Kolmo intro}) to a stochastic optimal control problem.
Let us consider the controlled equation
\begin{equation}
\left\{
\begin{array}
[c]{l}%
dX^{u}_\tau  =\left[  AX^{u}_\tau +\sqrt{Q} u_\tau  
 \right]  d\tau+\sqrt{Q}dW_\tau ,\text{ \ \ \ }\tau\in\left[  t,T\right] \\
X^{u}_t  =x.
\end{array}
\right.  \label{sdecontrolforte intro}%
\end{equation}
where the control $u$ takes values in a closed subset $K$ of $H$.
Define the cost
\begin{equation*}
J\left(  t,x,u\right)  =\mathbb{E}\int_{t}^{T}[
l\left(X^{u}_s\right)+g(u_s) ]ds+\mathbb{E}\phi\left(X^{u}_T\right). 
\end{equation*}
for real functions $l$, $\mathbb{\phi}$ and $g$ on $H$.
The control problem in strong formulation is to
minimize this functional $J$ over all admissible controls $u$.
We notice that we treat a control problem with unbounded controls,
and, in the case of superaquadratic hamiltonian, we require weak coercivity on the cost $J$.
Indeed, we assume that, for $1<q\leq 2$,
\[
0\leq g(u)\leq c(1+ \vert u \vert )^q, \quad\text{and}\quad g(u)\geq C \vert u \vert ^q \qquad \text{for every }u \in K : \vert u\vert \geq R
\]
so that the hamiltonian function 
\begin{equation*}
\psi\left(z\right)  =\inf_{u\in K}\left\{  g\left(u\right)
+zu\right\},\quad \forall z\in H ,
\end{equation*}
has quadratic growth in $z$ if $q=2$, and superquadratic growth of order $p>2$, the coniugate
exponent of $q$, if $q<2$.

Some example of operators $A$ and $Q$ commuting are listed in section \ref{sezionesde},
moreover this conditon is satisfied by a stochastic heat equation with coloured noise:
\begin{equation}\label{heat equation intro}
 \left\{
  \begin{array}{l}
  \dis
\frac{ \partial y}{\partial s}(s,\xi)= \Delta y(s,\xi)+
\frac{ \partial W^Q
}{\partial s}(s,\xi), \qquad s\in [t,T],\;
\xi\in \calo,
\\\dis
y(t,\xi)=x(\xi),
\\\dis
 y(s,\xi)=0, \quad \xi\in\partial \calo.
\end{array}
\right.
\end{equation}
Here $W^Q(s,\xi)$ is a Gaussian mean zero random field, such that
the operator $Q$ characterizes the correlation in the space variable.
The bounded linear operator $Q$ is diagonal with respect to
the basis $\left\lbrace e_k\right\rbrace_{k\in\N}$ of eigenfactors
of the Laplace operator with Dirichlet boundary conditions.
Equation (\ref{heat equation intro}) can be reformulated in $H=L^2(\calo)$
as an Ornstein-Uhlenbeck process (\ref{ornstein intro}) with $A$
and $Q$ commuting.

The paper is organized as follows:
in section \ref{sezionesde} some results on the Ornstein-Uhlenbeck process
are collected, in section \ref{sezioneforback} the Kolmogorov equation
(\ref{Kolmo intro}) is solved with $\psi$ with superquadratic growth
and $\phi$ lipschitz continuous, and these results are applied to
optimal control, in section \ref{sezioneforbackcont} the Kolmogorov equation
(\ref{Kolmo intro}) is solved with  quadratic $\psi$
and $\phi$ only continuous and again an application to control is briefly presented,
finally in section \ref{sez_contr_heat} optimal control problems for a
controlled heat equation are solved.
\section{Preliminary results on the forward equation and its semigroup}
\label{sezionesde}
We consider an Ornstein-Uhlenbeck process in a real and separable Hilbert space
$H$, that is a Markov process $X$ solution to equation%
\begin{equation}
\left\{
\begin{array}
[c]{l}%
dX_\tau  =AX_\tau d\tau+BdW_\tau
,\text{ \ \ \ }\tau\in\left[  t,T\right] \\
X_t =x,
\end{array}
\right.  \label{ornstein}%
\end{equation}
where $A$ is the generator of a strongly continuous semigroup in $H$ 
and $B$ is a linear bounded operator from $\Xi$ to $H$.
\noindent We define a positive and symmetric operator%
\[
Q_{\sigma}=\int_{0}^{\sigma}e^{sA}BB^{\ast}e^{sA^{\ast}}ds.
\]
Throughout the paper we assume the following.

\begin{hypothesis}
\label{ip su AB}

\begin{enumerate}
\item  The linear operator $A$ is the generator of a strongly continuous
semigroup $\left(  e^{t A},t\geq0\right)  $ in the Hilbert space $H.$
It is well known that there exist $M>0$ and $\omega\in\mathbb{R}$ such that
$\left\Vert e^{tA}\right\Vert _{L\left(  H,H\right)  }\leq Me^{\omega t}$, for
all $t\geq0$. In the following, we always consider $M\geq1$ and $\omega\geq 0$.
\item $B$ is a bounded linear operator from $\Xi$ to $H$ and $Q_{\sigma}$ is
of trace class for every $\sigma\geq0$.
\end{enumerate}
\end{hypothesis}
We notice that in some of the literature, in the case $\Xi=H,$ in order to define
the Ornstein-Uhlenbeck process, a bounded, symmetric and positive operator $Q$
is considered, and in equation \ref{ornstein}
$B$ is replaced by $\sqrt{Q}$.

\noindent The process $X$ is clearly time-homogeneous. It is well known that the
Ornstein-Uhlenbeck semigroup can be represented as $P_{\tau-t}=P_{t,\tau}$,
where%
\[
P_{\tau}\left[  \phi\right]  \left(  x\right)  :=\int_{H}\phi\left(  y\right)
\mathcal{N}\left(  e^{\tau A}x,Q_{\tau}\right)  \left(  dy\right)  ,
\]
and $\mathcal{N}\left(  e^{\tau A}x,Q_{\tau}\right)  \left(  dy\right)  $ denotes a
Gaussian measure with mean $e^{\tau A}x,$ and covariance operator $Q_{\tau}$.

In the following we are mainly concerned with the case $\Xi=H$, so we can take
$B=\sqrt{Q}$ and we assume that $A$ and $\sqrt{Q}$ commute. This happens
e.g. when $\left(  e_{n}\right)  _{n}$
is an orthonormal basis  in $H$ and $A$ and $Q$ have the spectral
decomposition $Ae_{n}=-\alpha_{n}e_{n}$ and $Qe_{n}=\gamma_{n}e_{n}$
where $\alpha_{n},\;\gamma_{n}%
>0$ and $\alpha_{n}\uparrow+\infty$. If the $\alpha_{n}$
are positive apart from a finite number the result is still true.

\noindent More in general let $\Xi=H$, $B=\sqrt{Q}$. Suppose that
$A$ is an unbounded, selfadjoint and negative defined operator,
$A=A^{\ast}\leq0,$ $A:\cald\left(  A\right)  \subset H\rightarrow
H,$ with inverse bounded. We consider the spectral representation of
$A$%
\begin{equation*}
A=\int_{-\infty}^{0}sdE\left(  s\right)  . 
\end{equation*}
and $Q=\left(  -A\right)  ^{\beta}=
{\displaystyle\int_{-\infty}^{0}}
\left(  -s\right)  ^{\beta}dE\left(  s\right)  $, for some $\beta
\in\mathbb{R}$. It turns out that
\[
Q_{t}=\int_{0}^{t}e^{sA}\left(  -A\right)  ^{\beta}e^{sA}ds=\frac{1}{2}\left(
1-e^{2At}\right)  \left(  -A\right)  ^{\beta-1}.
\]
It is proved in \cite{Mas} that in this case the Ornstein-Uhlenbeck semigroup
satisfies some regularizing property.  Namely the following hypothesis
\ref{ipH su fi} holds true with $\alpha=1/2  $.

\noindent We briefly introduce the notion of $\sqrt{Q}$-differentiability,
see e.g. \cite{Mas}. We recall that for a continuous function
$f:H\rightarrow\mathbb{R}$ the $\sqrt{Q}$-directional derivative $\nabla^{\sqrt{Q}}$ at a
point $x\in H$ in direction$\ \xi\in H$ is defined as follows:%
\[
\nabla^{\sqrt{Q}}f\left(  x;\xi\right)  =\lim_{s\rightarrow0}\frac{f\left(
x+s\sqrt{Q}\xi\right)  -f\left(  x\right)  }{s},\text{ }s\in\mathbb{R}\text{.}%
\]
A continuous function $f$ is $\sqrt{Q}$-G\^ateaux differentiable at a point $x\in H$ if
$f$ admits the $\sqrt{Q}$-directional derivative $\nabla^{\sqrt{Q}}f\left(  x;\xi\right)  $
in every directions $\xi\in H$ and there exists a functional, the
$\sqrt{Q}-$gradient $\nabla^{\sqrt{Q}}f\left(  x\right)  \in\Xi^{\ast}$ such that
$\nabla^{\sqrt{Q}}f\left(  x;\xi\right)  =\nabla^{\sqrt{Q}}f\left(  x\right)  \xi$. 
\begin{hypothesis}\label{ipH su fi} For some $\alpha\in [0,1)$ and for
every $\phi\in C_{b}\left(  H\right)  $, the function $P_{\tau}\left[
\phi\right]  \left(  x\right)  $ is $\sqrt{Q}  $-differentiable
with respect to $x$, for every $0\leq t <\tau < T$. Moreover there exists a
constant $c>0$ such that for every $\phi\in C_{b}\left(  H\right)  $, for
every $\xi\in\Xi,$ and for $0\leq t<\tau\leq T$,%
\begin{equation}
\left|  \nabla^{\sqrt{Q}}P_{\tau}\left[  \phi\right]
\left(  x\right)  \xi\right|  \leq\frac{c}{\tau  ^{\alpha}%
}\left\|  \phi\right\|  _{\infty}\left|  \xi\right|  .\label{ipotesi H}%
\end{equation}
\end{hypothesis}
In \cite{Mas} hypothesis \ref{ipH su fi} is verified by relating $\sqrt{Q}$-differentiability 
to properties of the operators $A$ and $\sqrt{Q}$. Namely if
\begin{equation}
\operatorname{Im}e^{tA}\sqrt{Q}\subset\operatorname{Im}Q_{t}^{1/2}.
\label{ornstein inclusione}%
\end{equation}
and for some $0\leq\alpha<1$ and $c>0$ the operator norm satisfies
\[
\left\|  Q_{t}^{-1/2}e^{tA}\sqrt{Q}\right\|  \leq ct^{-\alpha}\text{, for }0<t\leq
T.
\]
then hypotheses \ref{ipH su fi} is satisfied.

We also notice that this can be proved with a procedure similar to the one use in 
\cite{DP3} to prove the Bismut-Elworthy formula.
Namely, see e.g \cite{DP3}, lemma 7.7.2, for every uniformly continuous
function $\phi$ with bounded and uniformly continuous derivatives up to the second order,
\begin{equation}\label{Ito-mild}
\phi (X_\tau^{t,x})=P_{t,\tau}\phi(x)+\int_t^\tau<\nabla P_{s,\tau}\phi(X_s^{t,x}),\sqrt{Q}dW_s>
\qquad \P-\text{a.s.}
\end{equation}
By multiplying both sides of (\ref{Ito-mild}) by
\begin{equation*}
\int_t^\tau<\nabla (X_s^{t,x})h,dW_s>
\end{equation*}
and by taking expectation one gets
\begin{align}\label{pre-Bismut}
\E\left( \phi (X_\tau^{t,x})\int_t^\tau<\nabla (X_s^{t,x})h,dW_s>\right)
&=\E\left[ \int_t^\tau<\sqrt{Q}\nabla P_{s,\tau}\phi(X_s^{t,x}),\nabla (X_s^{t,x})h ds>\right] \\ \nonumber
&=E\left[ \int_t^\tau<\nabla P_{s,\tau}\phi(X_s^{t,x}),\sqrt{Q}e^{(s-t)A})h ds>\right] \\ \nonumber
&=\int_t^\tau<\nabla^{\sqrt{Q}} \E P_{s,\tau}\phi(X_s^{t,x}),h >ds \\ \nonumber
&=\int_t^\tau<\nabla^{\sqrt{Q}} P_{t,s} P_{s,\tau}\phi(x),h >ds \\ \nonumber
&=(\tau -t)<\nabla^{\sqrt{Q}} P_{t,\tau}\phi(x),h> \\ \nonumber
\end{align}
By arguments similar to the ones used in \cite{DP3}, lemma 7.7.5, we get that
for every bounded and continuous function $\phi$, the function $P_{\tau}\left[
\phi\right]  \left(  x\right)  $ is $\sqrt{Q}  $-differentiable
with respect to $x$, for every $0\leq t <\tau\leq T$ and 
\begin{equation*}
\left|  \nabla^{\sqrt{Q}}P_{t,\tau}\left[  \phi\right]
\left(  x\right)  \xi\right|  \leq\frac{c}{(\tau-t)  ^{1/2}%
}\left\|  \phi\right\|  _{\infty}\left|  \xi\right|  .
\end{equation*}
An analogous result can be proved in the case of $A$ and $Q$
commuting and $P$ transition semigroup of the perturbed Ornstein Uhlenbeck process
\begin{equation*}
\left\{
\begin{array}
[c]{l}%
dX_\tau  =AX_\tau d\tau+\sqrt{Q}F(\tau,X_\tau)+\sqrt{Q}dW_\tau
,\text{ \ \ \ }\tau\in\left[  t,T\right] \\
X_t =x,
\end{array}
\right.
\end{equation*}
with $F$ jointly continuous in $t$ and $x$ and lipschitz continuous in $x$ uniformly
with respect to $t$.

The model we have in mind consists of an heat equation.
Namely let $\calo$ be a bounded domain in $\R$.
We denote by $H$ the Hilbert space $L^2(\Omega)$ and by $\left\lbrace e_k\right\rbrace_{k\in\N}$
the complete orthonormal basis which diagonalizes $\Delta$,
endowed with Dirirchlet boundary conditions in $\calo$.
We consider the equation
\begin{equation}\label{heat equation}
 \left\{
  \begin{array}{l}
  \dis
\frac{ \partial y}{\partial s}(s,\xi)= \Delta y(s,\xi)+
\frac{ \partial W^Q
}{\partial s}(s,\xi), \qquad s\in [t,T],\;
\xi\in \calo,
\\\dis
y(t,\xi)=x(\xi),
\\\dis
 y(s,\xi)=0, \quad \xi\in\partial \calo.
\end{array}
\right.
\end{equation}
Here $W^Q(s,\xi)$ is a Gaussian mean zero random field, such that
the operator $Q$ characterizes the correlation in the space variables.
Namely the covariance of the noise is given by
$$
\E <W^Q(s, \cdot),h>_H <W^Q(t, \cdot),k>_H=t\wedge s<Qh,k>_H.
$$
\noindent In particular $W^Q(s,\xi)$ can be the Brownian sheet so that $\frac{ \partial^2 W^Q
}{\partial s \partial \xi}(s,\xi)$ in this case is the space-time white noise.
More in general we think about a coloured noise and on $Q$ we make the following assumptions:
\begin{hypothesis}\label{ip-cov}
 The bounded linear operator $Q:H\rightarrow H$ is positive and diagonal with respect to
the basis $\left\lbrace e_k\right\rbrace_{k\in\N}$, with eigenvalues
$\left\lbrace \lambda_k\right\rbrace_{k\in\N}$.
\end{hypothesis}
By previous assumptions it turns out that $\lambda_k\geq 0$.
Note that $W^Q(s,\cdot)$ is formally defined by
\[
 W^Q(s,\cdot)=\sum_{k=1}^n Qe_k(\cdot)\beta_k(s)
\]
where $\left\lbrace \beta_k(s)\right\rbrace_{k\in\N}$
is a sequence of mutually independent standard Brownian motions, all defined on the same
stochastic basis $(\Omega, \calf,\calf_t, \P)$.

Equation (\ref{heat equation}) can be written in an abstract way in $H$
as
\begin{equation}\label{heat eq abstr}
 \left\{
\begin{array}
[c]{l}%
dX_\tau  =AX_\tau d\tau+\sqrt{Q}dW_\tau
,\text{ \ \ \ }\tau\in\left[  t,T\right] \\
X_t =x,
\end{array}
\right.
\end{equation}
where $A$ is the Laplace operator with Dirirchlet boundary conditions,
$W$ is a cylindrical Wiener process in $H$ and $Q$ is its covariance operator.

\section{The semilinear Kolmogorov equation: lipschitz continuous final condition}
\label{sezioneforback}
The aim of this section is to present exitence and uniqueness results for the solution of
a semilinear Kolmogorov equation with the nonlinear term which is quadratic 
with respect to the $\sqrt{Q}$-derivative. The following arguments presented
in this section work also in the case
of $\psi$ with superquadratic growth with respect to $z$.

More precisely, let $\mathcal L $ be the generator of the transition
semigroup $P_t$, that is, at least formally, 
$$
(\call f)(x)=\frac{1}{2}(Tr Q^* \nabla^2 f)(x)+\<Ax,\nabla f(x)\>.
$$
Let us consider the following equation
\begin{equation}
\left\{
\begin{array}
[c]{l}%
\frac{\partial v}{\partial t}(t,x)=-\mathcal{L}v\left(  t,x\right)
+\psi\left(  \nabla^{\sqrt{Q}}v\left(  t,x\right)  \right)+l(x)  ,\text{ \ \ \ \ }t\in\left[  0,T\right]
,\text{ }x\in H\\
v(T,x)=\phi\left(  x\right)  ,
\end{array}
\right.  \label{Kolmo}%
\end{equation}
We introduce the notion of mild solution of the non linear Kolmogorov
equation (\ref{Kolmo}), see e.g. \cite{fute} and also \cite{Mas} for the definition of
mild solution when $\psi$ depends only on $\nabla^{\sqrt{Q}}v$ and not on $\nabla v$. 
Since $\mathcal{L}$ is (formally) the generator of
$P_{t}$, the variation of constants formula for (\ref{Kolmo}) is:%
\begin{equation}
v(t,x)=P_{t,T}\left[  \phi\right]  \left(  x\right)  -\int_{t}^{T}%
P_{t,s}\left[  \psi(\nabla^{\sqrt{Q} }v\left(  s,\cdot\right))  \right]  \left(  x\right)  ds,
-\int_{t}^{T}
P_{t,s}\left[  l\right]  \left(  x\right)  ds,\text{\ \ }t\in\left[  0,T\right]  ,\text{ }x\in H. \label{solmildkolmo}%
\end{equation}
We use this formula to give the notion of mild
solution for the non linear Kolmogorov equation (\ref{Kolmo}); we have also to
introduce some spaces of continuous functions, where we seek the
solution of (\ref{Kolmo}).

\noindent For $\alpha\geq0$, let $C_{\alpha}\left(  \left[  0,T\right]  \times
H\right)  $ (denoted by $C\left(  \left[  0,T\right]  \times
H\right)$ for $\alpha =0$) be the linear space of continuous functions $f:\left[  0,T\right)
\times H\rightarrow\mathbb{R}$ such that 
$$
\sup_{t\in\left[  0,T\right]  }%
\sup_{x\in H}\left(  T-t\right)  ^{\alpha}\left\vert f\left(  t,x\right)
\right\vert <+\infty.$$
$C_{\alpha}\left(  \left[  0,T\right]  \times H\right)$
endowed with the norm
\[
\left\Vert f\right\Vert _{C_{\alpha}}=\sup_{t\in\left[  0,T\right]  }%
\sup_{x\in H}\left(  T-t\right)  ^{\alpha}\left\vert f\left(  t,x\right)
\right\vert ,
\]
is a Banach space.

\noindent We consider also the linear space $C_{\alpha}^{s}\left(  \left[
0,T\right]  \times H,H^{\ast}\right)  $ (denoted by $C^s\left(  \left[  0,T\right]  \times
H,H^{\ast}\right)$ for $\alpha =0$) of the mappings $L:\left[
0,T\right)  \times H\rightarrow H^{\ast}$ such that for every $\xi\in H$,
$L\left(  \cdot,\cdot\right)  \xi\in C_{\alpha}\left(  \left[  0,T\right]
\times H\right)  $. The space $C_{\alpha}^{s}\left(  \left[  0,T\right]
\times H,H^{\ast}\right)  $ turns out to be a Banach space if it is endowed
with the norm
\[
\left\Vert L\right\Vert _{C_{\alpha}\left(  H^{\ast}\right)  }=\sup
_{t\in\left[  0,T\right]  }\sup_{x\in H}\left(  T-t\right)  ^{\alpha
}\left\Vert L\left(  t,x\right)  \right\Vert _{H^{\ast}}.
\]
In other words, $C_{\alpha}^{s}\left(  \left[  0,T\right]  \times H,H^{\ast
}\right)  $ can be identified with the space of the operators

\noindent$L\left( H,C_{\alpha}\left(  \left[  0,T\right]  \times H\right)
\right)  $.

\begin{definition}
\label{defsolmildkolmo}Let $\alpha\in\left[  0,1\right)  $. We say that a
function $v:\left[  0,T\right]  \times H\rightarrow\mathbb{R}$ is a mild
solution of the non linear Kolmogorov equation (\ref{Kolmo}) if the following
are satisfied:

\begin{enumerate}
\item $v\in C_{b}\left(  \left[  0,T\right]  \times H\right)  $;

\item $\nabla^{\sqrt{Q}  }v\in C_{\alpha}^{s}\left(  \left[
0,T\right]  \times H,H^{\ast}\right)  $: in particular this means that for
every $t\in\left[  0,T\right)  $, $v\left(  t,\cdot\right)$ is $\sqrt{Q}$-differentiable;

\item  equality (\ref{solmildkolmo}) holds.
\end{enumerate}
\end{definition}

Existence and uniqueness of a mild solution of equation (\ref{Kolmo})
is related to the study of the following 
forward-backward system: for given $t\in [0,T]$ and $x\in H$,
\begin{equation}\label{fbsde}
    \left\{\begin{array}{l}\dis dX_\tau =
AX_\tau d\tau+ \sqrt{Q} dW_\tau,\quad \tau\in
[t,T]\subset [0,T],
\\\dis
X_t=x,
\\\dis
 dY_\tau=-\psi(Z_\tau)\;d\tau-l(X_\tau)\;d\tau+Z_\tau\;dW_\tau,
  \\\dis
  Y_T=\phi(X_T),
\end{array}\right.
\end{equation}
and to the identification of $Z_t^{t,x}=\nabla_x Y_t^{t,x} \sqrt{Q} $.
We extend the definition of $X$ setting
$X_s=x$ for $0\leq s\leq t$. The second equation in
(\ref{fbsde}), namely
\begin{equation}\label{bsde}
    \left\{\begin{array}{l}\dis
 dY_\tau=-\psi(Z_\tau)\;d\tau-l(X_\tau)\;d\tau+Z_\tau\;dW_\tau,
 \qquad \tau\in [0,T],
  \\\dis
  Y_T=\phi(X_T),
\end{array}\right.
\end{equation}
is of backward type.
Under suitable assumptions on the coefficients
 $\psi:H
\rightarrow\mathbb{R}$, $l:H\rightarrow\mathbb{R}$
and $\mathbb{\phi}:H\rightarrow\mathbb{R}$
we will look for a solution consisting of a pair of predictable processes,
taking values in $\mathbb{R}\times H$, such that $Y$ has
continuous paths and
\[
\|\left( Y,Z\right)\|^2_{\mathbb{K}_{cont}}:=
\mathbb{E}\sup_{\tau\in\left[ 0,T\right] }\left\vert Y_{\tau}\right\vert
^{2}+\mathbb{E}\int_{0}^{T}\left\vert Z_{\tau}\right\vert ^{2}d\tau<\infty,
\]
see e.g. \cite{PaPe1}. In the following we denote by $\mathbb{K}_{cont}\left(
\left[ 0,T\right] \right)$ the space of such processes.

The solution of (\ref{fbsde}) will be denoted by $(X_\tau, Y_\tau, Z_\tau)_{\tau\in[0,T]}$, or,
to stress the dependence on the initial time $t$ and on the
initial datum $x$, by $(X_\tau^{t,x}, Y_\tau^{t,x}, Z_\tau^{t,x})_{\tau\in[0,T]}$.
In the following we refer to \cite{fute} for the definition of the class 
$\calg(H)$ of G\^ateaux differentiable functions $f:H\rightarrow \R$
with strongly continuous derivative.
\begin{hypothesis}
\label{ip su psi}The maps $\psi:H
\rightarrow\mathbb{R}$, $l:H\rightarrow\mathbb{R}$ and $\mathbb{\phi}:H\rightarrow\mathbb{R}$ are Borel
measurable, moreover $\psi$ is G\^ateaux differentiable, namely 
$\psi\in \calg(H)$ () and for every $\xi_1,\xi_2\in \Xi$,
$\vert \psi(\xi_1)-\psi(\xi_2)\vert\leq (1+\vert\xi_1\vert^{p-1}
+\vert\xi_2\vert^{p-1})\vert\xi_1-\xi_2\vert $, for $p\geq 2$. The maps 
$l$ and $\phi$ belong to $C_b(H)$.

\noindent Moreover from now on we assume that, unless modifying the value of
$l$, $\psi(0)=0$.
\end{hypothesis}

We make differentiability assumptions on the coefficients of equation 
(\ref{bsde}):
\begin{hypothesis}\label{ipfidiffle}
 $l$ and $\phi$ belong to $\calg(H)$ and have bounded derivative.
\end{hypothesis}

Assume that in hypothesis \ref{ip su psi} $p=2$.
So by \cite{Kob}, under hypothesis \ref{ip su psi} the BSDE (\ref{bsde}) admit a unique solution and by
\cite{BriFu}, under the further assumption \ref{ipfidiffle} setting $v(t,x):=Y_t^{t,x}$,
it turns out that $v$ is the unique mild solution of equation (\ref{Kolmo}), 
and $\nabla^{\sqrt{Q}}v(t,x)=Z_t^{t,x}.$
By assuming that $A$ and $Q$ commute, or, more in general,
by assuming that hypothesis \ref{ipH su fi} holds true, also imposing a more restrive structure on the
forward equation and on the backward equation, we will prove in section \ref{sezioneforbackcont} an estimate on $Z_\tau^{t,x}$
depending on $\tau$, $t$, $T$ and $\Vert \phi \Vert_\infty$ but not on $\nabla \phi$.
Thanks to this estimate we will prove that by setting 
\begin{equation}\label{def-v}
v(t,x):=Y_t^{t,x},
\end{equation}
it turns out that $v$ is the unique mild solution of equation (\ref{Kolmo}), 
and $\nabla^{\sqrt{Q}}v(t,x)=Z_t^{t,x}$ without assumption \ref{ipfidiffle}.
We note that differentiability on $l$, thanks to the regularizing property of
the semigroup, can be easily removed. So from now on we can consider the case of $l=0$.

We go on in this section with the study of equation (\ref{Kolmo}) with $\psi$ superquadratic.
\subsection{Local mild solution of the corresponding PDE}
\label{subsezPDE}

In this subsection we look for a local mild solution of equation (\ref{Kolmo}), that
is a mild solution in a small time interval.
We work on the PDE following \cite{Go2},
but the same result can be achieved by working on the BSDE with a procedure similar to the one indicate in \cite{BaoDeHu}, section 4.1.

We start by proving existence of a local mild solution to equation (\ref{Kolmo}),
and then we look for a priori estimates for this local mild solution.

\begin{theorem}\label{teokolmo}
Assume that hypotheses \ref{ip su AB}, \ref{ipH su fi}, \ref{ip su psi}, \ref{ipfidiffle}
hold true. Then equation (\ref{Kolmo}) admits a unique local mild solution $u$
on $[T-\delta, T]$, for some $0<\delta<T$ according to definition \ref{defsolmildkolmo}. 
\end{theorem}

\noindent {\bf Proof.} The first part of the proof is similar to
the proof of theorem 2.9 in
\cite{Mas} and we partially omit it. Consider the product space 
$\Lambda:C\left(  \left[  0,T\right]  \times H\right)  \times
C^s\left(  \left[  0,T\right]  \times H,H^{\ast}\right)  $ with
the product norm. Let us also denote by $\Lambda_{R_0}$
the closed ball of radius $R_0$, with respect to the product norm 
\begin{equation*}
 \|(f,L)\|_{C\left(  \left[  0,T\right]  \times H\right)  \times
C^s\left(  \left[  0,T\right]  \times H,H^{\ast}\right)}=
\|f\|_{C\left(  \left[  0,T\right]  \times H\right) }+
\|L\|_{C^s\left(  \left[  0,T\right]  \times H,H^{\ast}\right)}.
\end{equation*}

% which is 
% the product of the equivalent norms in $C\left(  \left[  0,T\right]  \times H\right)$
% and in $C^s\left(  \left[  0,T\right]  \times H,\Xi^{\ast}\right)$:
% 
% \begin{align*}
%   &\left\|  f\right\|  _{\beta,C }=\sup_{t\in\left[  0,T\right]  }%
% \sup_{x\in H}\exp\left(  -\beta\left(  T-t\right)  \right) \left| f(t,x)  \right| ,\\ \nonumber
% & \left\|  L\right\|  _{\beta,C^{s}\left(  \Xi^{\ast}\right)  }%
% =\sup_{t\in\left[  0,T\right]  }\sup_{x\in H}\exp\left(  -\beta\left(
% T-t\right)  \right)  \left\|  L\left(
% x,t\right)  \right\|  _{\Xi^{\ast}},
% \end{align*}
% where $\beta$ is a positive constant to be fixed.

Let us also define, for $(u,v)\in \Lambda_{R_0}$ 
\begin{equation}
\Gamma_{1}\left[  v,w\right]  (t,x)=P_{t,T}\left[  \phi\right]  \left(
x\right)  -\int_{t}^{T}P_{t,s}\left[  \psi\left( w\left( 
s,\cdot\right)  \right)  \right]  \left(  x\right)ds
-\int_{t}^{T}P_{t,s}\left[  l  \right]  \left(  x\right)ds
, \label{gamma1}%
\end{equation}%
\begin{equation}
\Gamma_{2}\left[  v,w\right]  (t,x)=\nabla^{\sqrt{Q} }%
P_{t,T}\left[ \phi\right]  \left(  x\right)  -\int_{t}^{T}\nabla^{\sqrt{Q}}P_{t,s}\left[  
\psi\left( w\left(  s,\cdot\right)  \right)  \right]  \left(  x\right)  ds.
-\int_{t}^{T}\nabla^{\sqrt{Q}} P_{t,s}\left[  l  \right]  \left(  x\right)ds
\label{gamma2}%
\end{equation}
Thanks to condition (\ref{ipH su fi}) and if $\delta$ is sufficiently small,
$\Gamma$ is well defined on $\Lambda_{R_0}$ with
values in itself.
Indeed, let us take 
\[
 R_0=2Me^{\omega T}(\|\phi\|_1+T\|l\|_1)
\]
It turns out that
\begin{align*}
 &\|\Gamma_{1}\left[  v,w\right] , \Gamma_{2}\left[  v,w\right]\|
_{C\left(  \left[  0,T\right]  \times H\right)  \times
C^s\left(  \left[  0,T\right]  \times H,\Xi^{\ast}\right)} \\ \nonumber
&\leq (\|\phi\|_1+T\|l\|_1+(T-t)C_{R_0}R_0+(T-t)^{1-\alpha}CC_{R_0})R_0 \\ \nonumber
&\leq (\frac{1}{2}+(T-t)C_{R_0}+(T-t)^{1-\alpha}CC_{R_0})R_0, \\ \nonumber
\end{align*}
where $C_{R_0}$ is the lipschitz constant of $\psi$ if
$\vert v\vert_{C^s\left(  \left[  0,T\right]  \times H,H^{\ast}\right)}\leq R_0$.
For $t\in[T-\delta,T]$ and $\delta$ sufficiently small
we get that $\frac{1}{2}+\delta C_{R_0}+\delta^{1-\alpha}CC_{R_0}<1$,
so that $\Gamma:\Lambda_{R_0}\rightarrow\Lambda_{R_0}$
 Moreover, arguing as in \cite{Mas}, theorem 2.9, 
it is possible to show that $\Gamma$ is a contraction in 
$\Lambda_{R_0}$, 
and so we are able to find a unique local mild solution to equation
(\ref{Kolmo}). \qed

\subsection{Equivalent representation of the mild solution}
\label{subsezequiv}
In this subsection we give an an alternative representation of the mild solution
of equation (\ref{Kolmo}).
Let $v$ be the local
mild solution of equation (\ref{Kolmo}), as stated in theorem \ref{teokolmo}.
Let us define
\begin{equation}\label{defG}
G(t,x)=\int_0^1 \nabla \psi(\lambda\nabla^{\sqrt{Q}}v(t,x))d\lambda. 
\end{equation}
We present in an informal way the object of this subsection. Equation
\eqref{Kolmo} can be rewritten as
\begin{equation}
\left\{
\begin{array}
[c]{l}%
\frac{\partial v}{\partial t}(t,x)=-\mathcal{L}v\left(  t,x\right)
+\< \sqrt{Q}G(t,x),\nabla v(t,x) \> +l(x)  ,\text{ \ \ \ \ }t\in\left[  0,T\right]
,\text{ }x\in H\\
u(T,x)=\phi\left(  x\right)  ,
\end{array}
\right.  \label{Kolmobis}%
\end{equation}
Let us consider the Markov process $\Theta^{t,x}_\cdot$ solution to equation 
\begin{equation}
\left\{
\begin{array}
[c]{l}%
d\Theta_\tau^{t,x} =A\Theta_\tau^{t,x} d\tau 
+\sqrt{Q}G(\tau,\Theta_\tau^{t,x})d\tau+\sqrt{Q}dW_\tau
,\text{ \ \ \ }\tau\in\left[  t,T\right]  \text{ \ \ \ }  t\in[T-\delta,T]\\
\Theta_t =x,
\end{array}
\right.  \label{ornsteinpert}%
\end{equation}
Let us denote by $R_{t,\tau}$ the transition semigroup of $\Theta^{t,x}_\cdot$
Following \cite{Go2}, since the operator $\mathcal{L}u\left(  t,x\right)
+\< \sqrt{Q}G(t,x),\nabla u(t,x) \>$ is formally the generator of $R$, 
we want to prove that the mild solution of
(\ref{Kolmo}) for $t\in[T-\delta,t]$ can be represented as
\begin{equation}
v(t,x)=R_{t,T}\left[  \phi\right]  \left(  x\right) 
-\int_{t}^{T}
R_{t,s}\left[  l\right]  \left(  x\right)  ds,\text{\ \ }t\in\left[T-\delta,T\right]  ,\text{ }x\in H. \label{solmildkolmobis}%
\end{equation}
Representation (\ref{solmildkolmobis}) immediately gives an a priori
estimate for the norm of $v$ in $C([0,T],H)$.

In the following lemma we notice that equation (\ref{ornsteinpert})
admits a unique mild solution in weak sense.
First of all we state an existence and uniqueness result for equation (\ref{ornsteinpert}).
\begin{lemma}\label{lemmaornsteinpert}
Let hypothesis \ref{ip su AB} and \ref{ip su psi} on 
$\psi$ hold true, let $G$ be defined by (\ref{defG}), then equation (\ref{ornsteinpert})
has a unique mild solution in weak sense, and this solution is unique in law.
\end{lemma}
 \dim
By the Girsanov theorem, since $\vert G(t,x) \vert\leq R_0$, there exists a probability
measure 
$\tilde{\P}$, equivalent to the original one $\P$, such that
\[
 \left\lbrace \tilde{W}_\tau=\int_0^\tau G(r,x) dr+W_\tau, \tau\geq 0 \right\rbrace 
\]
is a Brownian motion.
In the probability space $(\Omega, \calf, \tilde{\P})$, $\Theta^{t,x}$
is an Ornstein-Uhlenbeck process:
\begin{equation*}
\left\{
\begin{array}
[c]{l}%
d\Theta^{t,x}_\tau  =A\Theta^{t,x}_\tau d\tau+\sqrt{Q}d\tilde{W}_\tau
,\text{ \ \ \ }\tau\in\left[  t,T\right] \\
\Theta^{t,x}_t =x,
\end{array}
\right.  %
\end{equation*}
and this guarantees existence and uniqueness in law of a weak solution to
equation \ref{ornsteinpert}.
This suffices to have the transition semigroup $R$ well defined.
\qed

Next we want to prove that representation (\ref{solmildkolmobis}) holds true.
A similar result is obtained in \cite{Go2} by using the results in \cite{Ce}, \cite{CeGo} about Cauchy problems 
associated to weakly continuous semigroups, such as transition semigroup.
Here we use the connection between PDEs and BSDEs. We notice that
in \cite{BaoDeHu} a BSDE in a Markovian framework with generator
with superqadratic growth is solved also in the case of final datum continuous in $x$.
With the following techniques we solve a semilinear Komogorov equation like (\ref{Kolmo}),
in the case of lipschitz continuous final datum. Notice that by asking some smoothing properties
of the transition semigroups allows also in this case to identify in the corresponding BSDE
$Z_t^{t,x}$ with $\nabla ^{\sqrt{Q}}Y_t^{t,x}$. In this way our results can be applied
to solve a related stochastic optimal control problem.
\begin{proposition}\label{proprappreq}
 Assume that hypotheses \ref{ip su AB}, \ref{ipH su fi}, 
\ref{ip su psi}, \ref{ipfidiffle} hold true.
Then the local mild solution $v$ of (\ref{Kolmo}) can be represented as in
(\ref{solmildkolmobis}).
\end{proposition}
\dim
Let $u$ be the local mild solution of (\ref{Kolmo}) and let us define
%for $0\leq t \leq \tau \leq T$
\begin{equation}\label{solbsde}
Y_\tau^{t,x}=v(\tau,X_\tau^{t,x}), \text{\ \ \ }
Z_\tau^{t,x}=\nabla^{\sqrt{Q}} v(\tau,X_\tau^{t,x}),
\end{equation}
where as usual for $0\leq \tau\leq t$, $X_\tau^{t,x}=x$.
It is well known that the pair of processes $(Y_\tau^{t,x},Z_\tau^{t,x})_{0\leq \tau \leq T}$
is solution of the BSDE (\ref{bsde}). Moreover since $v$ is a local mild solution 
of (\ref{Kolmo}), $\vert Z_\tau^{t,x} \vert \leq R_0$ $\P$-a.s.. Since by our assumption \ref{ip su psi}
$\psi $ is locally lipschitz continuous it turns out that
\[
 f(\tau):=
\left\lbrace
\begin{array}[l]{ll}
 \dfrac{\psi(Z^{t,x}_\tau)}{\vert Z^{t,x}_\tau \vert^2}Z^{t,x}_\tau & \text{ if }Z^{t,x}_\tau\neq0 \\
0 & \text{ otherwise }
\end{array}
\right. 
\]
is bounded. So, by using the techniques introduced in \cite{BriHu},
by the Girsanov theorem there exists a probability measure
$\tilde{\P}$, equivalent to the original one $\P$, such that
\[
 \left\lbrace \tilde{W}_\tau=-\int_0^\tau f(r) dr+W_\tau, \tau\geq 0 \right\rbrace 
\]
is a Brownian motion.
So in $(\Omega, \calf, \tilde{\P})$ equation (\ref{bsde}) can be rewritten as
\begin{equation}\label{bsdetilde}
    \left\{\begin{array}{l}\dis
 dY_\tau=-l(X_\tau)\;d\tau+Z_\tau\;d\tilde{W}_\tau,
 \qquad \tau\in [0,T],
  \\\dis
  Y_T=\phi(X_T).
\end{array}\right.
\end{equation}
We notice that in $(\Omega, \calf, \tilde{\P})$ $X$ is solution to
\begin{equation}
\left\{
\begin{array}
[c]{l}%
dX_\tau  =AX_\tau d\tau +\sqrt{Q}G(\tau,X_\tau)d\tau+\sqrt{Q}dW_\tau
,\text{ \ \ \ }\tau\in\left[  t,T\right] \\
X_t =x,
\end{array}
\right.  \label{ornsteintilde}%
\end{equation}
and lemma \ref{lemmaornsteinpert} guarantees existence, uniqueness
and regularity of the mild solution of this equation. Moreover with respect to the new probability
measure $\tilde{\P}$ the transition semigroup  of $X$ coincides with $R_{t,T}$. In particular we
notice that the law of $(X,Y,Z)$ depends on the coefficients $A$, $G$, $\psi$, $l$, $\phi$, on the initial
condition $x$ given at initial time $t$, but not on the probability space nor on the Wiener
process. So in particular it turns out that again, since $Y^{t,x}_t$ is deterministic,
it coincides with $u(t,x)$, mild solution of equation
\ref{Kolmobis}. So it turns out that mild representation \ref{solmildkolmobis}
holds true.
\qed

\begin{remark}\label{remark_lip}
It is possible to prove that for the local mild solution
$u$ of equation (\ref{Kolmo}) there exists $L>0$
such that for every $x,y\in H$ and for every $t\in[T-\delta, T]$
\begin{equation}
(T-t)^\alpha \vert \nabla^{\sqrt{Q}} v(t,x)-\nabla^{\sqrt{Q}} v(t,y)\vert \leq L \vert x-y\vert
\label{ulipschitz}
\end{equation}
This can be proved by showing that $\Gamma$ defined in the
proof of theorem \ref{teokolmo} is a contraction
in the product space $C\left(  \left[  0,T\right]  \times H\right)  \times
C^s_{1,\alpha}\left(  \left[  0,T\right]  \times H,H^{\ast}\right)  $
where by $C^s_{1,\alpha}\left(  \left[  0,T\right]  \times H,H^{\ast}\right)$
we mean the space of the operators $L\in C^{s}_\alpha\left(  \left[  0,T\right]  \times H,H^{\ast}\right)$
such that $(T-t)^\alpha L(t,x) $ is lipschitz continuous in $x$
uniformly with respect to $t$.
We endow $C^s_{1,\alpha}\left(  \left[  0,T\right]  \times H,H^{\ast}\right)$
with the norm
$$
\left\Vert L\right\Vert _{C_{1,\alpha}\left( H^{\ast}\right)  }=\sup
_{t\in\left[  0,T\right]  }\sup_{x\in H}
\left\Vert L\left(  t,x\right)  \right\Vert _{H^{\ast}}
+\sup_{t\in\left[  0,T\right]  }\sup_{x,y\in H}\left(  T-t\right)  ^{\alpha
}\left\Vert L\left(  t,x\right)- L\left(  t,y\right) \right\Vert _{H^{\ast}}.
$$
and we endow $C\left(  \left[  0,T\right]  \times H\right)  \times
C^s_{1,\alpha}\left(  \left[  0,T\right]  \times H,H^{\ast}\right)  $
with the product norm.
 Let us also denote by $\Lambda_{R_0}$
the closed ball of radius $R_0$, with respect to the product norm 
\begin{equation*}
 \|(f,L)\|_{C\left(  \left[  0,T\right]  \times H\right)  \times
C_{1,\alpha}\left(  \left[  0,T\right]  \times H,H^{\ast}\right)}=
\|f\|_{C\left(  \left[  0,T\right]  \times H\right) }+
\|L\|_{C_{1,\alpha}\left(  \left[  0,T\right]  \times H,H^{\ast}\right)}.
\end{equation*}
Let us also define $\Gamma_{1}$ and $\Gamma_{2}$
as in (\ref{gamma1}) and (\ref{gamma2}).
We take $R_0$ and $\delta$ as in the proof of theorem \ref{teokolmo}.
We have to prove that $\Gamma:\Lambda_{R_0}\rightarrow\Lambda_{R_0}$
and it is a contraction. For the firts point, in view of the results of theorem
\ref{teokolmo} we have to prove that $(T-t)^\alpha\Gamma_{2}\left[  v,w\right] (t,\cdot)$
is lipschitz continuous.
 For $\xi\in H$ we have
\begin{align*}
 &\vert (T-t)^\alpha(\Gamma_{2}\left[  v,w\right](t,x)-
\Gamma_{2}\left[  v,w\right](t,y))\vert \\ \nonumber
&\vert(T-t)^\alpha \nabla^{\sqrt{Q}} P_{t,T}[\phi](x)\xi-\nabla^{\sqrt{Q}} P_{t,T}[\phi](y)\xi \vert
+(T-t)^\alpha\int_t^T\vert \nabla^{\sqrt{Q}} P_{t,s}[l](x)\xi-\nabla^{\sqrt{Q}} P_{t,s}[l](y)\xi \vert ds \\ \nonumber
&+(T-t)^\alpha\int_t^T\vert \nabla^{\sqrt{Q}} P_{t,s}[\psi(w(s,\cdot))](x)\xi
-\nabla^{\sqrt{Q}} P_{t,s}[w(s,\cdot))](y) \xi\vert ds
=I+II+III. \\ \nonumber
\end{align*}

We start by estimating $I$, following e.g. \cite{DP1} and \cite{Mas} we get :
\begin{align*}
I & 
=(T-t)^\alpha\vert\nabla^{\sqrt{Q}} \int_H ( \phi(z+e^{(T-t)A}x)-\phi(z+e^{(T-t)A}y) )
\caln (0,Q_{T-t}) dz \vert \\ \nonumber
&=(T-t)^\alpha\vert \int_H ( \phi(z+e^{(T-t)A}x)-\phi(z+e^{(T-t)A}y) )
\< Q_{T-t}^{-1/2}e^{(T-t)A}\sqrt{Q}\xi, Q_{T-t}^{-1/2}z \>\caln (0,Q_{T-t}) (dz)  \vert \\ \nonumber
&\leq (T-t)^\alpha Me^{\omega (T-t)}\|  \phi\| _1 \vert x-y \vert _H \|  Q_{T-t}^{-1/2}e^{(T-t)A}\sqrt{Q}\|
\vert \xi \vert\leq Me^{\omega (T-t)}\Vert \phi\Vert_1 \vert x-y \vert _H 
\vert \xi \\ \nonumber
\end{align*}
Arguing in a similar way it follows that 
$$
II\leq (T-t)Me^{\omega (T-t)}\Vert \phi\Vert_1
$$ 
so that 
$$
I+II\leq R_0/2
$$
For what concerns $III$ we get
\begin{align*}
III 
& = (T-t)^\alpha\int_t^T\int_H \vert \psi(w(s,z+e^{(s-t)A}x))-
 \psi(w(s,z+e^{(s-t)A}y))\vert \\ \nonumber
&\< Q_{s-t}^{-1/2}e^{(s-t)A}\sqrt{Q}\xi, Q_{s-t}^{-1/2}z \>\caln (0,Q_{s-t}) (dz)ds \\ \nonumber
&\leq (T-t)^\alpha R_0 C_{R_0} Me^{\omega(T-t)} \vert x-y \vert _H 
\int_t^T (T-s)^{-\alpha}\|Q_{s-t}^{-1/2}e^{(s-t)A}\sqrt{Q}\|
\vert \xi \vert ds \\ \nonumber
&\leq 2R_0 C_{R_0} (T-t)^{1-\alpha}, \\ \nonumber
\end{align*}
where we have used the fact that $u$ is the local mild solution and $\psi$
is locally lipschitz continuous.
Also by the proof of theorem \ref{teokolmo} it turns out that
\begin{align*}
 &\|\Gamma_{1}\left[  v,w\right] , \Gamma_{2}\left[  v,w\right]\|
_{C\left(  \left[  0,T\right]  \times H\right)  \times
C^{s,1,\alpha}\left(  \left[  0,T\right]  \times H,H^{\ast}\right)} \\ \nonumber
&\leq (\|\phi\|_1+T\|l\|_1+(T-t)C_{R_0}R_0+(T-t)^{1-\alpha}CC_{R_0})R_0 \\ \nonumber
&\leq (1+(T-t)C_{R_0}+3(T-t)^{1-\alpha}CC_{R_0})R_0, \\ \nonumber
\end{align*}
Let $0<\bar\delta\leq\delta$ be such that for $t\in[T-\bar\delta,T]$ 
we get that $1+\bar\delta C_{R_0}+3\bar\delta^{1-\alpha}CC_{R_0}<1$,
so that $\Gamma:\Lambda_{R_0}\rightarrow\Lambda_{R_0}$.

Moreover, it is possible to show that $\Gamma$ is a contraction in 
$\Lambda_{R_0}$, 
and so we are able to find a unique local mild solution to equation
(\ref{Kolmo}) in $C\left(  \left[  0,T\right]  \times H\right)  \times
C^{s,1,\alpha}\left(  \left[  0,T\right]  \times H,H^{\ast}\right)$.

\qed
\end{remark}

As a consequence equation (\ref{ornsteinpert}) admits a
unique mild solution, in classical (strong) sense.
\begin{lemma}\label{lemmaornsteinpert_forte}
Let hypothesis \ref{ip su AB} and \ref{ip su psi} on 
$\psi$ hold true, let $G$ be defined by (\ref{defG}), then equation (\ref{ornsteinpert})
has a unique mild solution satisfying moreover, for every $x,y \in H$,
$$
\vert \Theta^{t,x}_\tau -\Theta^{t,y}_\tau  \vert \leq C_T
\vert x-y  \vert 
$$
\end{lemma}
\dim
The proof is standard apart from the singularity of the Lipschitz constant of $G$ in $T=t$:
by theorem \ref{teokolmo}, the local mild solution $u$ is Lipschitz continuous according to
estimate \ref{ulipschitz}, and so also $G$ is. Indeed, for every $x,y \in H$,
\begin{equation*}
 (T-t)^\alpha \vert G(t,x)-G(t,y)\vert \leq L \vert x-y\vert
\label{Glipschitz}
\end{equation*}
Existence, uniqueness and Lipschitz property of a mild solution 
of equation \ref{ornsteinpert} follow as in \cite{Go2}, proposition 3.9.
\qed

The aim of the next section is to find a priori estimates for the local 
mild solution of equation \ref{Kolmo} by using reperesentation (\ref{solmildkolmobis}).
We notice that, the transition semigroup $R_{t,T}$
is a perturbed Ornstein-Uhlenbeck transition semigroup, so
we could try to investigate if it satisfies regularizing properties 
like the ones satisfied by the Ornstein-Uhlenbeck transition
semigroup contained in \ref{ipH su fi}.
Anyway there are some difficulties related to the coefficient
$G$ in equation (\ref{ornsteinpert}): $G$ is
not differentiable and blows up 
like $(T-t)^{-\alpha}$, so it is in general not
square integrable since $0<\alpha<1$. In particular when $A$ and $Q$ commute,
hypothesis \ref{ipH su fi} holds true with $\alpha=1/2$. In the existing literature, see e.g.
\cite{fu}, \cite{Mas1} and references therein,
regularizing properties for perturbed 
Ornstein-Uhlenbeck transition semigroup,
such as the strong Feller property or property 
\ref{ipH su fi}, are proved by means of ``generalizations"
of the Girsanov theorem and then by means of the Malliavin calculus, eventually
with direct calculation of the Malliavin derivative. Here this cannot be done:
since $G$ is not square integrable, no immediate generalization
of the Girsanov theorem can be applied. Also, $G$ is not differentiable,
so the existing theory cannot be directly used, even if
in this direction generalizations seem less involved.

\subsection{A priori estimates and global existence}
\label{sezapriori}

In this section we investigate a priori estimates for the local 
mild solution of equation \ref{Kolmo} that we have found in theorem \ref{teokolmo}. 
By proposition \ref{proprappreq}, the equivalent representation (\ref{solmildkolmobis})
for the mild solution $v$ holds true and this immediately gives an a priori estimate for 
the supremum norm of $v$, namely
\begin{equation}\label{apriori_cont}
\left\|  v\right\|  _{C }\leq  (\left\|  \phi \right\|_0+ T \left\|  l\right\| _0)\leq \frac{R_0}{2}.
\end{equation}
Next, we look for a priori estimates for the norm
$\left\| \cdot\right\|  _{C\left( H^{\ast}\right)  }$ of $\nabla^G v$.

\noindent In order to prove an a priori estimate for the norm 
$\left\| \cdot\right\|  _{C\left( H^{\ast}\right)  }$ of $\nabla^G v$,
we exploit again the strict connection between PDEs and BSDEs.
We start by the fact the if $v$ is the local mild solution of (\ref{Kolmo}),
and $X_\cdot^{t,x}$ is the Ornstein-Uhlenbeck process defined 
by (\ref{ornstein}), then $((v(\tau,X_\tau^{t,x} ))_\tau,(\nabla v(\tau,X_\tau^{t,x}))_\tau)$
is solution to the the BSDE (\ref{bsde}), that we rewrite for the reader convenience:
\begin{equation}\label{bsde3}
    \left\{\begin{array}{l}\dis
 dY^{t,x}_\tau=-\psi(Z^{t,x}_\tau)\;d\tau-l(X^{t,x}_\tau)\;d\tau+Z^{t,x}_\tau\;dW_\tau,
 \qquad \tau\in [0,T],
  \\\dis
  Y^{t,x}_T=\phi(X^{t,x}_T),
\end{array}\right.
\end{equation}
 For $\xi\in H$ let us define, if it exists, $F^{t,x}_\tau=\nabla^{\sqrt{Q}} Y^{t,x}_\tau \xi$,
$V^{t,x}_\tau=\nabla^{\sqrt{Q}} Z^{t,x}_\tau \xi$. It turns out that such processes exist,
$(F^{t,x}_\tau , V^{t,x}_\tau)\in \mathbb{K}_{cont}\left(
\left[ 0,T\right] \right)$, and that they are solution to 
\begin{equation}\label{bsde_der}
    \left\{\begin{array}{l}\dis
 dF^{t,x}_\tau=-\nabla\psi(Z^{t,x}_\tau)V^{t,x}_\tau\;d\tau
-\nabla l(X^{t,x}_\tau)e^{(\tau-t)A}{\sqrt{Q}}\xi\;d\tau+V^{t,x}_\tau\;dW_\tau,
 \qquad \tau\in [0,T],
  \\\dis
  F^{t,x}_T=\nabla \phi(X^{t,x}_T)e^{(T-t)A}{\sqrt{Q}}\xi ,
\end{array}\right.
\end{equation} which is equation \ref{bsde3}
differentiated in direction $\sqrt{Q}\xi$, since for $t\leq\tau\leq T$,
$\nabla^{\sqrt{Q}} X^{t,x}_\tau\xi=e^{(\tau-t)A}\sqrt{Q}\xi$.
\begin{proposition}\label{propbsde_der}
Assume that hypotheses \ref{ip su AB}, \ref{ipH su fi}, \ref{ip su psi}
and \ref{ipfidiffle} hold true. Then equation (\ref{bsde_der}) admits
a unique solution that is a pair of predictable processes,
taking values in $\mathbb{R}\times H$, such that $F$ has
continuous paths and
% \[
% \|\left( F, V\right)\|_{\mathbb{K}_{cont}}<\infty.
% \]
% In particular 
$F_t^{t,x}$ is bounded.
\end{proposition}
\dim
We notice that again $Z^{t,x}_\tau=\nabla^{\sqrt{Q}} u(\tau, X^{t,x}_\tau)$
where $u$ is the local mild solution of equation \ref{Kolmo}. Since
$\nabla\psi$ is locally lipschitz continuous, it turns out that
\[
 f_1(\tau):=
\left\lbrace
\begin{array}[l]{ll}
 \dfrac{\nabla\psi(Z^{t,x}_\tau)}{\vert Z^{t,x}_\tau \vert^2}Z^{t,x}_\tau & \text{ if }Z^{t,x}_\tau\neq0 \\
0 & \text{ otherwise },
\end{array}
\right. 
\]
is bounded. Following again the method
in \cite{BriHu} by the Girsanov theorem there exists a probability measure
$\hat{\P}$, equivalent to the original one $\P$, such that
\[
 \left\lbrace \hat{W}_\tau=-\int_0^\tau f_1(r) dr+W_\tau, \text{ \ \ }
\tau\geq 0 \right\rbrace 
\]
is a Brownian motion.
So in $(\Omega, \calf, \hat{\P})$ equation (\ref{bsde_der}) can be rewritten as
\begin{equation}\label{bsde_der1}
    \left\{\begin{array}{l}\dis
 dF^{t,x}_\tau=-\nabla l(X^{t,x}_\tau)e^{(\tau-t)A}{\sqrt{Q}}\xi\;d\tau
+Z^{t,x}_\tau\;d\hat{W}_\tau,
 \qquad \tau\in [0,T],
  \\\dis
  F^{t,x}_T=\nabla\phi(X_T^{t,x})e^{(T-t)A}{\sqrt{Q}}\xi.
\end{array}\right.
\end{equation}
In this equation the generator $-\nabla l(X^{t,x}_\tau)e^{(\tau-t)A}{\sqrt{Q}}\xi$
is independent on $F$ and $V$ and it is bounded, so by classical theorems on
BSDEs equation (\ref{bsde_der1}) admits a unique solution $(F,V)$
such that $F$ has continuous paths and
\[
\|\left( F,V\right)\|^2_{\hat{\mathbb{K}}_{cont}}:=
\hat{\mathbb{E}}\sup_{\tau\in\left[ 0,T\right] }\left\vert F^{t,x}_{\tau}\right\vert^{2}
+\hat{\mathbb{E}}\int_{0}^{T}\left\vert V^{t,x}_{\tau}\right\vert ^{2}d\tau<C,
\]
see e.g. \cite{PaPe1}. Notice that the constant $C$ depends only on $A$, $G$,
$l$, $\psi$ on the initial condition $x$ given at initial time $t$. In particular 
\[
\vert F^{t,x}_t \vert < Me^{\omega T}(\|\phi\|_1+T\|l\|_1)
\]
\qed

We immediately deduce the following result:
\begin{corollary}\label{cor_aprioriest}
Assume that hypotheses \ref{ip su AB}, \ref{ipH su fi}
%(or equivalently \ref{ipABequivalenti})
, \ref{ip su psi},
\ref{ipfidiffle} hold true and let $u$ be the local mild solution 
of equation (\ref{Kolmo}), as stated in theorem \ref{teokolmo}. 
\begin{equation}\label{apriori}
\left\|  v\right\|  _{C([0,T]\times H) }+\left\| \nabla^{\sqrt{Q}} v\right\|  _{C\left([0,T]\times H, H^{\ast}\right)  }
\leq R_0.
\end{equation}

\end{corollary}
\dim
The estimate for the norm $\left\|  \cdot\right\|  _{C }$ of $u$
follows by (\ref{apriori_cont}), the estimate for the norm
$\left\| \cdot\right\|  _{C\left( H^{\ast}\right)  }$
of $\nabla^{\sqrt{Q}}u$ follows by proposition \ref{propbsde_der}.
\qed

We can state a result on existence and uniqueness of a mild solution $u$ of
equation (\ref{Kolmo}), which immediately gives a unique mild solution
of equation (\ref{bsde}).
\begin{theorem}\label{teobsde}
Assume that hypotheses \ref{ip su AB}, \ref{ipH su fi}, \ref{ip su psi} and \ref{ipfidiffle}
hold true. Then equation (\ref{Kolmo}) admits a unique mild solution $u$
according to definition \ref{defsolmildkolmo}. 
Let $X_\cdot^{t,x}$ be solution of equation (\ref{ornstein}).
The pair of processes $(Y_\tau^{t,x}=v(\tau,X_\tau^{t,x}),
Z_\tau^{t,x}=\nabla^{\sqrt{Q}} v(\tau,X_\tau^{t,x}))_{\tau\in[0,T]}$ is the unique solution 
of the BSDE (\ref{bsde}).
\end{theorem}
\dim
The global existence of the mild solution $v$ follows by the local existence
(Theorem \ref{teokolmo}) and by the a priori estimates (Corollary \ref{cor_aprioriest}).
The connections between PDEs and BSDEs is classical in the literature
also in the infinite dimensional case
(see e.g. \cite{fute}) and the proof is complete.
\qed

The next step is to remove differentiability assumptions on $l$
and $\phi$. We start by assuming $l$ bounded and continuous
and $\phi$ bounded and lipschitz continuous.

\begin{theorem}\label{teobsdelip}
Assume that hypotheses \ref{ip su AB}, \ref{ipH su fi}, \ref{ip su psi},
hold true and that $l$ is bounded and continous and $\phi$
is bounded and lipschitz continuous.
Then equation (\ref{bsde}) admits a unique solution,
that is a pair of processes  $(Y_\cdot^{t,x},
Z_\cdot^{t,x})\in \mathbb K_{cont}([0,T])$.
The function $v(t,x)=Y_t^{t,x}$
is the unique mild solution of equation (\ref{Kolmo})
according to definition \ref{defsolmildkolmo}.
\end{theorem}
\dim
We consider the inf-sup convolution of $\phi$ (see e.g. \cite{LL} 
and \cite{DP3})
denoted by $\phi_{n}$ and defined by
\begin{equation}
\phi_{n}\left(  x\right)  =\sup_{z\in
H}\left\{  \inf_{y\in H}\left[  \phi\left(  y\right)  +n\frac{\left|
z-y\right|  _{H}^{2}}{2}\right]  -n\left|  x-z\right|  _{H}^{2}
\right\}  . \label{infsupconv}%
\end{equation}
Similarly, let us define $l_{n}$ the inf-sup convolution of $l$.
It is well known that $\phi_n\in UC_{b}^{1,1}\left(  H\right)$
and as $n$ tends to $+\infty$, $\phi
_{n}$ converges to $\phi$ uniformly. 
Moreover, see also \cite{Mas}, let us denote by $L$
the Lipschitz constant of $\phi$; then $\left| \nabla \phi\right|\leq L$.

\noindent Now let us denote by $(Y_\cdot^{n,t,x},
Z_\cdot^{n,t,x})$ the unique solution
 of the BSDE (\ref{bsde}) with $\phi_n$ and $l_n$ in the
place of $\phi$ and $l$ respectively.
By standard results on BSDEs we know that as $n\rightarrow \infty$
\[
\E \sup_{\tau\in\left[ 0,T\right] }\left\vert Y_\tau^{n,t,x}-Y_\tau^{t,x}\right\vert
^{2}+\mathbb{E}\int_{0}^{T}\left\vert Z_\tau^{n,t,x}-Z_\tau^{t,x}\right\vert ^{2}d\tau
\rightarrow 0
\]
and moreover 
\[
\|\left( Y^n,Z^n\right)\|^2_{\mathbb{K}_{cont}}< C
\]
where $C$ is a constant independent on $n$.
\noindent We need to prove some further regularity on $Z$
let us denote by $(F_\cdot^{n,t,x},
V_\cdot^{n,t,x})$ the unique solution
 of the BSDE (\ref{bsde_der}) with $\phi_n$ and $l_n$ in the
place of $\phi$ and $l$ respectively.
It turns out that 
\[
\E \sup_{\tau\in\left[ 0,T\right] }\left\vert F_\tau^{n,t,x}\right\vert^2<C
\]
where $C$ is a constant depending on $L$ and on 
$\Vert\phi \Vert_\infty$, and independent on $n$.
So the process $Z_\cdot^{n,t,x}$ is uniformly bounded in 
$n$, since for $\xi \in H$, $Z_\tau^{n,t,x}\xi= F_\tau^{n,t,x}$.
We also get that
\[
\E \sup_{\tau\in\left[ 0,T\right] }\left\vert Y_\tau^{t,x}\right\vert^2+
\E \sup_{\tau\in\left[ 0,T\right] }\left\vert Z_\tau^{t,x}\right\vert^2<C.
\]
By setting $v(t,x)=Y_t^{t,x}$, $\nabla^Gv(t,x)=Z_t^{t,x}$,
we have found a (unique ) mild solution to (\ref{Kolmo}).
\qed

\subsection{Application to control}
\label{applic contr 1}

We formulate the stochastic optimal control problem in the strong
sense. Let $\left(  \Omega,\mathcal{F},\mathbb{P}\right)  $ be a given
complete probability space with a filtration $\left(  \mathcal{F}_{\tau
}\right)  _{\tau\geq0}$ satisfying the usual conditions. $\left\{  W\left(
\tau\right)  ,\tau\geq0\right\}  $ is a cylindrical Wiener process on $H
$\ with respect to $\left(  \mathcal{F}_{\tau}\right)  _{\tau\geq0}$. The
control $u$ is an $\left(  \mathcal{F}_{\tau}\right)  _{\tau}$-predictable
process with values in a closed set $K$ of a normed space $U$; in the following
we will make further assumptions on the control processes.
Let $R: U\rightarrow H$
and consider the controlled state equation
\begin{equation}
\left\{
\begin{array}
[c]{l}%
dX^{u}_\tau  =\left[  AX^{u}_\tau +\sqrt{Q}R\left( u_\tau)  \right)
 \right]  d\tau+\sqrt{Q}dW_\tau ,\text{ \ \ \ }\tau\in\left[  t,T\right] \\
X^{u}_t  =x.
\end{array}
\right.  \label{sdecontrolforte}%
\end{equation}
The solution of this equation will be denoted by
$X_\tau^{u,t,x}$ or simply by $X^{u}_\tau$. $X$ is also called
the state, $u$ and $T>0,$ $t\in\left[  0,T\right]$ are fixed.
The special structure of equation (\ref{sdecontrolforte}) allows to study
the optimal control problem related by means of BSDEs and
(\ref{sdecontrolforte}) leads to a semilinear Hamilton Jacobi Bellman
equation which is a special case of the Kolmogorov equation (\ref{Kolmo}) we
have studied in the previous sections. The occurrence of the operator
$\sqrt{Q}$ in the control term is imposed by our techniques, on the contrary the
presence of the operator $R$ allows more generality.

Beside equation (\ref{sdecontrolforte}),\ define the cost
\begin{equation}
J\left(  t,x,u\right)  =\mathbb{E}\int_{t}^{T}[
l\left(X^{u}_s\right)+g(u_s) ]ds+\mathbb{E}\phi\left(X^{u}_T\right). 
\label{cost}%
\end{equation}
for real functions $l$, $\mathbb{\phi}$ on $H$ and $g$ on $U$.
The control problem in strong formulation is to
minimize this functional $J$ over all admissible controls $u$.
We make the following assumptions on the cost $J$.

\begin{hypothesis}
\label{ip costo}

\begin{enumerate}
\item  The function $\mathbb{\phi}:H\rightarrow\mathbb{R}$ is lipschitz
continuous and bounded;

\item $l: H\rightarrow\mathbb{R}$ is bounded and continuous;

\item $g: U\rightarrow\mathbb{R}$ is mesurable; and for some $1<q\leq 2$ there exists a constant $c>0$
such that 
\begin{equation}
 \label{crescita costo1}
0\leq g(u)\leq c(1+ \vert u \vert ^q) 
\end{equation}
and there exist $R>0$, $C>0$ such that
\begin{equation}
\label{crescita costo}
g(u)\geq C \vert u \vert ^q \qquad \text{for every }u \in K : \vert u\vert \geq R.
\end{equation}

\end{enumerate}
\end{hypothesis}
In the following we denote by 
$\mathcal{A}_{d}$ the set of
admissible controls, that is the $K$-valued predictable processes such that 
\[
 \E \int_0^T \vert u_t \vert ^q dt <+\infty.
\]
This summability requirement is justified  by (\ref{crescita costo}):
a control process which is not $q$-summable would have infinite cost.

\noindent We denote by $J^{\ast}\left(  t,x\right)  =\inf_{u\in\mathcal{A}%
_{d}}J\left(  t,x,u\right)  $ the value function of the problem and, if it
exists, by $u^{\ast}$ the control realizing the infimum, which is called
optimal control.

\noindent We make the following assumptions on $R$.

\begin{hypothesis}
\label{ip aggiuntive}
$R: U\rightarrow H$ is measurable and $\vert R(u)\vert\leq C(1+\vert u\vert )$
for every $u\in U$.
\end{hypothesis}

We have to show that equation (\ref{sdecontrolforte}) admits a unique mild solution, for
every admissible control $u$.
\begin{proposition} \label{prop_sde_contr}
 Let $u$ be an admissible control and assume that hypothesis \ref{ip su AB} holds true.
Then equation (\ref{sdecontrolforte}) amits a unique mild solution $(X_\tau)_{\tau\in[t,T]}$
such that $\E\sup_{\tau\in[t,T]}\vert X_\tau\vert^q< \infty$.
\end{proposition}

{\bf Proof.} The proof follows the proof of proposition 2.3 in
\cite{fuhute}, with suitable changes since in that paper
the finite dimensional case in considered and the current cost $g$
has quadratic growth with respect to $u$, that is in (\ref{crescita costo}) $q=2$. 

In order to make an approximation procedure in (\ref{sdecontrolforte})
we introduce the sequence of stopping times
$$
\tau_n=\inf \left\lbrace t\in[0,T]: \E \int_0^t \vert u_s \vert ^q ds>n \right\rbrace 
$$
with the ususal convention that $\tau_n=T$ if this set is empty. Following the
approximation procedure used in the proof of proposition 2.3 in
\cite{fuhute} we can prove that there exists a unique mild solution
with the required $q$-integrability.
\qed

We define in a classical way the Hamiltonian function relative to the above
problem:%
\begin{equation}\label{hamilton}
\psi\left(z\right)  =\inf_{u\in K}\left\{  g\left(u\right)
+zR(u)\right\}\quad \forall z\in H .
\end{equation}
We prove that the hamiltonian function just defined satisfies the polynomial growth
conditions and the local lipschitzianity required in hypothesis \ref{ip su psi}.
\begin{lemma}
 \label{lemma-hamilton}
The hamiltonian $\psi:H\rightarrow \R$ is Borel measurable, there exists a constant $C>0$ such that
\[
-C(1+\vert z\vert^p)\leq \psi(z)  \leq g(u)+\vert z\vert (1+\vert u\vert), \qquad \forall u \in K,
\]
where $p$ is the coniugate exponent of $q$.
Moreover if the infimum in (\ref{hamilton}) is attained, it is attained
in a ball of radius $C(1+\vert z\vert^{p-1})$ that is 
\[
\psi(z)  =\inf_{u\in K,\vert u\vert \leq C(1+\vert z\vert^{p-1})}\left\{  g\left(u\right)
+zR(u)\right\},\quad z \in H,
\]
and
\[
\psi(z)  <  g\left(u\right)
+zR(u)\quad \text{if }\vert u\vert > C(1+\vert z\vert^{p-1}).
\]
In particular it follows that $\psi$ is locally lipschitz continuous, namely
$\forall\,z_1,z_2\in H$, for some $C>0$,
\begin{equation}\label{stima2psi}
 \vert \psi(z_1)-\psi(z_2) \vert\leq C(1+\vert z_1\vert^{p-1}+\vert z_2 \vert^{p-1})\vert z_1-z_2 \vert
\end{equation}
\end{lemma}
{\bf Proof.}
The measurability of $\psi$ is straightforward. By assumption (\ref{crescita costo})
we get 
\begin{equation}\label{stima1psi}
g(u)+zR(u)\geq C(\vert u\vert^q-\vert R\vert^q)-C_1\vert z\vert(1+\vert u\vert)
\end{equation}
where $C$ and $R$ are as in (\ref{crescita costo}) and $C_1>0$,
and by this it follows that 
$$
\psi\left(z\right)  \geq\inf_{u\in U}\left\{  g\left(u\right)
+zR(u)\right\}\geq C(\vert u\vert^q-\vert R\vert^q)-C_1\vert z\vert(1+\vert u\vert)
\geq -C_2 \vert z\vert^p -C_3
$$
for suitable constants $C_2$ and $C_3$.
Moreover
$$
\vert\psi(z)\vert\leq g(u)+c\vert z\vert(1+\vert u\vert).
$$
Now we prove that the infimum is attained in the ball of radius
$C(1+\vert z\vert^{p-1})$. By (\ref{stima1psi}),
\begin{equation*}
 g(u)+zR(u)\geq C\vert u \vert(\vert u\vert^{q-1}-\frac{C_1}{C}\vert z\vert)
-C\vert R\vert^q- C_1\vert z\vert.
\end{equation*}
On the other hand, for some $u^0\in K:$
$$
g(u^0)+zR(u^0)\leq C_4(1+\vert z\vert).
$$
and so there exists a constant $\bar C$ such that
if $\vert u \vert \geq \bar C(1+\vert z\vert^{p-1})$ then
$$
g(u)+zR(u)\geq g(u^0)+zR(u^0)
$$
and the result follows from the continuity of $g$ and $R$.
Finally (\ref{stima2psi}) now easily follows.
\qed

\begin{remark}
 We give an example of hamiltonian we can treat.
Let $g(u)=\vert u\vert ^q$, $1<q\leq 2$. Then, if $R(u)=u$, the hamiltonian function turns out to be
\[
 \psi(z)=\left( \left(\dfrac{1}{q}\right)^{1/(q-1)}-\left(\dfrac{1}{q}\right)^p\right)\vert z\vert^p 
\]
where $p\geq2$ is the coniugate of $q$.
With this example, for $p=2$, the Hamilton Jacobi Bellman related
can be solved with the ad hoc exponential transform, see e.g.\cite{Go2}.
Our theory cover also the case of
hamiltonian functions not exactly equal to $\vert z \vert ^2$.
\end{remark}
We define
\begin{equation}\label{defdigammagrande}
\Gamma(s,x,z)=\left\{ u\in U: zR(u)+g(u)= \Psi(z)\right\};
\end{equation}
if $\Gamma(z) \neq \emptyset$ for every $z\in H$, by \cite{AuFr}, see Theorems 8.2.10 and
8.2.11, $\Gamma$ admits a measurable selection, i.e. there exists
a measurable function $\gamma: H \rightarrow U$ with
$\gamma(z)\in \Gamma(z)$ for every $z\in \R$.

In the following theorem, in order to prove the so called fundamental relation,
we have to make further assumptions concerning differentiability
of the hamiltonian function $\psi$.
These assumptions allow us to say that the Hamilton Jacobi Bellman equation
relative to the above problem, which is given by equation (\ref{Kolmo}),
admits a unique mild solution by theorem \ref{teokolmo}.
Moreover this solution can be represented by means
of the solution of the BSDE
(\ref{bsde}), namely the solution is given by
$v(t,x)=Y_t^{t,x}$. So, adequating to our context
the techniques e.g. in \cite{fuhute},
we can prove the fundamental relation for the optimal control.

\begin{theorem}\label{th-rel-font}
 Assume hypotheses \ref{ip su AB}, \ref{ipH su fi}, 
\ref{ip costo} and \ref{ip aggiuntive} hold true, and assume that the hamiltonian function $\psi$
satisfies G\^ateaux differentiability assumptions stated in hypothesis \ref{ip su psi}.
 For
every $t\in [0,T]$, $x\in H$ and for
 all admissible control $u$ we have $J(t,x,u(\cdot))
 \geq v(t,x)$,
  and the
 equality holds if and only if
$$
u_s\in \Gamma\left( \nabla^{\sqrt{Q}}
v(s ,X^{u,t,x}_s)
\right)
  $$
\end{theorem}
{\bf Proof.}
% Now we prove a fundamental relation when the hamiltonian function has
% superquadratic growth, that is $p>2$ in lemma \ref{lemma-hamilton} (and $1<q<2$
% for what concerns the growth of the cost, see hypothesis \ref{ip costo}). We notice that the following arguments also work in the quadratic case.
% 
For every admissible control $(u_t)_{t\in[0,T]}$, we define, for every $n\in \N$, 
\[
u^n_t=u_t 1_{\vert u_t\vert\leq n}+n 1_{\vert u_t\vert >n}.
\]
Since $u\in L^q(\Omega\times [0,T])$, then $u^n\rightarrow u$ in
$L^q(\Omega\times [0,T])$ and so $u^n\rightarrow u$ with respect to the measure $dt\times\P$.
Since moreover the sequence $(u_n)_n$ is monotone, then the convergence holds
for almost all $t\in[0,T]$ and $\P$-almost surely.
Moreover it follows that 
\[
\int_0^T \vert R( u^n_s) \vert^2 ds
\leq C\int_0^T(1+ \vert u^n_s \vert)^2 ds
\leq C(1+n^2).
\]
Let us define 
\[
\rho_n=\exp\left( - \int_0^T R( u^n_s) dW_s-\frac{1}{2}\int_0^T \vert R( u^n_s) \vert^2 ds \right) 
\]
In the probability space $(\Omega, \calf, \P)$ let $X^{u^n}$ denote the solution
of equation 
\begin{equation}\label{sdecontrol_n}
 \left\{
\begin{array}
[c]{l}%
dX^{u^n}_\tau  =\left[  AX^{u^n}_\tau +\sqrt{Q}R\left( u^n_\tau)  \right)
 \right]  d\tau+\sqrt{Q}dW_\tau ,\text{ \ \ \ }\tau\in\left[  t,T\right] \\
X^{u^n}_t  =x.
\end{array}
\right. 
\end{equation}
By the Girsanov theorem there exists a probability measure $\P^n$,
equivalent to the original one $\P$, namely $\dfrac{d\P^n}{d\P}=\rho_n$,
and such that $W^n_t:=W_t+\int_0^t R(u^n_s) ds$
is a $\P^n$-Wiener process. In $(\Omega, \calf, \P^n)$, $X^{u^n}$ solves the
following stochastic differential equation
\begin{equation}\label{sdecontrol_n_Gir}
 \left\{
\begin{array}
[c]{l}%
dX^{u^n}_\tau  = AX^{u^n}_\tau +\sqrt{Q}dW^n_\tau ,\text{ \ \ \ }\tau\in\left[  t,T \right] \\
X^{u^n}_t  =x.
\end{array}
\right. 
\end{equation}
Let us also denote by $(Y^n,Z^n)$ the solution in $(\Omega, \calf, \P^n)$ of the BSDE
\begin{equation}\label{bsde-control}
    \left\{\begin{array}{l}\dis
 dY^n_\tau=-\psi(Z^n_\tau)\;d\tau-l(X^{u^n}_\tau)\;d\tau+Z^n_\tau\;dW^n_\tau,
 \qquad \tau\in [0,T],
  \\\dis
  Y_T=\phi(X^{u^n}_T).
\end{array}\right.
\end{equation}
We notice that the law of $(X^{u^n},Y^n,Z^n)$ depends on the coefficients $A$, $\sqrt{Q}$, $\psi$,
$l$ and $\phi$ and not on the Wiener process nor on the probability space.
So in particular $Y^{n,t,x}_t=v(t,x) $, where $v$ has been defined
in (\ref{def-v}) and it is also the solution of the related Hamilton-Jacobi-Bellman equation.
This fact will be crucial in order to study the convergence of $(Y^n,Z^n)$
as $n\rightarrow +\infty$.

\noindent In the probability space $(\Omega, \calf, \P)$, we denote by $X^u$
the solution of equation (\ref{sdecontrolforte}), then 
\begin{align*}
 \E\sup_{t\in[0,T]}\vert X^{u^n}-X^u\vert^q &\leq \E\sup_{t\in[0,T]}\vert\int_0^te^{(t-s)A}\sqrt{Q}
(u_s-n)1_{\vert u_s \vert >n}ds\\ \nonumber
&\leq C(T,A,Q)\E\int_0^T \vert u_s-n\vert^q 1_{\vert u_s \vert>n}ds\\ \nonumber
\end{align*}
and so $X^{u^n}\rightarrow X^u$ in $L^q(\Omega,C([0,T],H))$, with probability measure $\P$.
By combining this fact with
the previous arguments on the the $dt\times\P$-
almost sure convergence of $u^n$ we get that $X^{u^n}_t\rightarrow X^u_t$ $\P$-almost surely uniformly with respect to $t$.

\noindent Moreover $(Y^{n,t,x}_t,Z^{n,t,x}_t)=(v(t,x),\nabla^{\sqrt{Q}}v(t,x)$,
and for all $\tau \in [t,T]$,
\begin{equation}
\label{identific_n}
 (Y^{n,t,x}_\tau,Z^{n,t,\tau}_t)=(v(\tau,X^{u^n,t,x}_\tau),\nabla^{\sqrt{Q}}v(\tau,X^{u^n,t,x}_\tau)
\end{equation}
and since we work with lipschitz continuous assumptions on the final cost
$\phi$ and G\^ateaux differentiability assumptions
on $\psi$ by theorem \ref{teokolmo} we get that both $v$ and $\nabla^{\sqrt{Q}}v$ are bounded and continuous.
In the probability space $(\Omega, \calf, \P)$, by
using the pointwise convergence, of $X^{u^n}$ to $X^u$,
uniformly with respect to time, we deduce a pointwise convergence uniform in time of
$(Y^n,Z^n)$, passing through the identification (\ref{identific_n}) of $Y^n$ and $Z^n$.
Namely, in $(\Omega, \calf, \P)$, $(Y^{n,t,x}_\tau,Z^{n,t,x}_\tau)=
(v(\tau,X^{u^n,t,x}_\tau),\nabla^{\sqrt{Q}}v(\tau,X^{u^n,t,x}_\tau))\rightarrow 
(v(\tau,X^{t,x}_\tau),\nabla^{\sqrt{Q}}v(\tau,X^{t,x}_\tau))$
$\P$-almost surely uniformly with respect to $\tau$.

Now we are ready to prove the fundamental relation:
we integrate the BSDE (\ref{bsde-control}) in $[t,T]$: at first we write down the equation
with rspect to the $\P^n$-Wiener process $W^n$ and then we 
pass to the process $W$, which is a standard Wiener process
in the original probability space $(\Omega, \calf, \P)$:
\begin{align*}
 dY^n_t &=\phi(X^{u^n}_T)+\int_t^T\psi (Z^n_s)\;ds+
\int_t^T l(X^{u^n}_s)\;ds-\int_t^T Z^n_s\;dW^n_s \\ \nonumber
& =\phi(X^{u^n}_T)+\int_t^T\psi (Z^n_s)\;ds+\int_t^T l(X^{u^n}_s)\;ds
-\int_t^T R(u^n_s) ds -\int_t^T Z^n_s\;dW_s  \nonumber
\end{align*}
We notice that by standard arguments since
$Z^n\in L^2((\Omega,\calf,\P^n)\times[0,T])$, then it also holds that
$Z^n\in L^2((\Omega,\calf,\P)\times[0,T])$. Now we integrate with respect to
the original probability $\P$: by taking expectation in the previous integral equality
\begin{align*}
 dY^n_t &=\E\phi(X^{u^n}_T)+\E\int_t^T\psi (Z^{n}_s)\;ds
+\E\int_t^T l(X^{u^n}s)\;ds-\E\int_t^T Z^n_s\;dW^n_s \\ \nonumber
&= \E\phi(X^{u^n}_T)+\E\int_t^T\psi (Z^n_s)\;ds+\E\int_t^Tl(X^{u^n}_s)\;ds
-\E\int_t^T Z^n_s R(u^n_s) ds -\E\int_t^T Z^n_s\;dW_s \\ \nonumber
&= \E\phi(X^{u^n}_T)+\E\int_t^T\psi (\nabla^{\sqrt{Q}}v(s,X^{u^n}_s))\;ds+\E\int_t^T l(X^{u^n}_s)\;ds
-\E\int_t^T \nabla^{\sqrt{Q}}v(s,X^{u^n}_s) R(u^n_s) ds  \\ \nonumber
\end{align*}
where in the last passage the stochastic integral has zero expectation,
and we have identified $Z^n_s$ with $\nabla^{\sqrt{Q}}v(s,X^{u^n}_s)$.
Next we also identify $Y^n_t$ with $v(t,x)$ and then we let $n\rightarrow +\infty$:
\begin{align*}
 v(t,x) &
= \E\phi(X^{u^n}_T)+\E\int_t^T\psi (\nabla^{\sqrt{Q}}v(s,X^{u^n}_s))\;ds+\E\int_t^Tl(X^{u^n}_s)\;ds
-\E\int_t^T \nabla^{\sqrt{Q}}v(s,X^{u^n}_s) R(u^n_s) ds  \\ \nonumber
&\rightarrow \E\phi(X^u_T)+\E\int_t^T\psi (\nabla^{\sqrt{Q}}v(s,X^u_s))\;ds+\E\int_t^Tl(X^u_s)\;ds
-\E\int_t^T \nabla^{\sqrt{Q}}v(s,X^u_s) R(u_s) ds. \\ \nonumber
\end{align*}
By adding and subtracting $\E\int_t^T g(u_s) ds$ we get
\begin{equation}\label{rel-fond-dim}
 J(t,x,u) =v(t,x)+\E\int_t^T\left[ -\psi (\nabla^{\sqrt{Q}}v(s,X^u_s))
+\nabla^{\sqrt{Q}}v(s,X^u_s) R(u_s) +g(u_s)\right] ds,
\end{equation}
from which we deduce the desired conclusion.
 \qed

Under the assumptions of Theorem
\ref{th-rel-font}, let us define the so called
optimal feedback law:
\begin{equation}\label{leggecontrolloottima}
u(s,x)=\gamma\Big(\nabla^{\sqrt{Q}}
v(s ,X^{u,t,x}_s) \Big),\qquad
s\in [t,T],\;x\in H.
\end{equation}
Assume that the closed loop equation admits a solution
$\{\overline{X}_s,\;s\in
[t,T]\}$:
\begin{equation}\label{cle}
\overline{X}_s= e^{(s-t)A}x_0+
\int_t^s e^{(s-r)A}\sqrt{Q}dW_r
+\int_{t}^s e^{(s-r)A}R(\gamma(\nabla^{\sqrt{Q}}
v(s ,\overline{X}_r)))dr,
\qquad s\in [t,T].
\end{equation}
Then the pair $(\overline{u}=u(s,\overline{X}_s),\overline{X})_{s\in[t,T]}$
is optimal for the coNntrol problem.
We nevertheless notice that existence of a solution of the closed loop
equation is not obvious, due to the lack
of regularity of the feedback law $u$ occurring in
(\ref{cle}).
This problem can be avoided by formulating the optimal control problem in the weak sense, following \cite{FlSo}, see also \cite{fute} and \cite{Mas}.

By an {\em admissible control system}
we mean
$$(\Omega,\mathcal{F},
\left(\mathcal{F}_{t}\right) _{t\geq 0}, \mathbb{P}, W,
u(\cdot),X^u),$$
where $W$ is an $H$-valued Wiener process, $u$ is an admissible control and $X^u$
solves the controlled equation (\ref{sdecontrolforte}).
The control problem in weak formulation is to minimize
the cost functional over all the admissible control systems.

\begin{theorem}\label{teo su controllo debole}
Assume hypotheses \ref{ip su AB}, \ref{ipH su fi},
\ref{ip costo} and \ref{ip aggiuntive} hold true,
and assume that the hamiltonian function $\psi$
satisfies G\^ateaux differentiability assumptions
stated in hypothesis \ref{ip su psi}.
 For
every $t\in [0,T]$, $x\in H$ and for
 all admissible control systems we have $J(t,x,u(\cdot))
 \geq v(t,x)$,
  and the
 equality holds if and only if
$$
u_s\in \Gamma\left( \nabla^{\sqrt{Q}}
v(s ,X^{u}_s)
\right)
  $$

Moreover
assume that the set-valued map $\Gamma$ is non empty and let $\gamma$
be its measurable selection. 
\begin{equation*}
u_{\tau}=\gamma(\nabla^{\sqrt{Q}} v(\tau,X^u_\tau ))
,\text{ \ }\mathbb{P}\text{-a.s. for a.a. }\tau
\in\left[ t,T\right]
\end{equation*}
is optimal.

Finally, the closed loop equation 
\begin{equation}\label{closed loop eq}
 \left\{
\begin{array}
[c]{l}%
dX^{u}_\tau  =\left[  AX^{u}_\tau +\sqrt{Q
}R\left( \gamma(\nabla^{\sqrt{Q}} v(\tau,X^u_\tau )) )  \right)
 \right]  d\tau+\sqrt{Q}dW_\tau ,\text{ \ \ \ }\tau\in\left[  t,T\right] \\
X^{u}_t  =x.
\end{array}
\right.
\end{equation}
admits a weak solution
$(\Omega,\mathcal{F},
\left(\mathcal{F}_{t}\right) _{t\geq 0}, \mathbb{P}, W,X)$
which is unique in law and setting
$$
u_{\tau}=\gamma\left(\nabla^{\sqrt{Q}} v(\tau,X_\tau )\right),
$$
we obtain an optimal admissible control system $\left(
W,u,X\right) $.
\end{theorem}

\noindent {\bf Proof.} The proof follows from the fundamental relation
stated in theorem \ref{th-rel-font}.
The only difference here is the solvability
of the closed loop equation in the weak sense: this is a standard application of the 
Girsanov theorem. Indeed, by lemma \ref{lemma-hamilton}, the infimum in the hamiltonian
is achieved in a ball of radius $C(1+\vert z \vert ^{p-1})$ and so for the optimal control
$u$ the following estimate holds true, $\P$-a.s. and for a.a. $\tau\in[t,T],\;0\leq t\leq T$:
\[
 \vert u_\tau \vert \leq C(1+\vert Z_\tau^{t,x}\vert^{p-1})
=C(1+\vert \nabla^{\sqrt{Q}}v(\tau,X_\tau^{t,x}\vert^{p-1})\leq \bar C.
\]
Thanks to this bound we can apply a Girsanov change of measure
and the conclusion follows in a standard way.
\qed

\section{The semilinear Komogorov equation in the quadratic case: continuous final condition}
\label{sezioneforbackcont}

Let $v$ be the mild solution of the semilinear Kolmogororv equation \ref{Kolmo},
with the nonlinear term which is quadratic  with respect to the $\sqrt{Q}$-derivative
and with final condition $\phi$ differentiable.
The aim of this section is to present an estimate for the $\sqrt{Q}$-derivative
of $u(t,x)$ depending on $T,\;t,\;\Vert\phi\Vert_\infty$ but not
on the $\nabla \phi$. If $X$ is finite dimensional, also for processes
more general than the Onstein uhlenbeck process, this estimate has been obtained in 
$\cite{Ri}$ by imposing some conditions on the coefficient of the forward equation
 for $X$. With those conditions, by inverting $\nabla X$ and by techniques coming from BMO
martingales, the estimate is proved. Also in \cite{BaoDeHu} a similar estimate is proved
with $X$ finite dimensional, with a more restrictive structure than in \cite{Ri} but with
$\psi$ also with superquadratic growth. Applications of this estimate
to a related Kolmogorov equation are not exploited in \cite{BaoDeHu} nor in \cite{Ri}.

Here we prove the estimate for $\nabla^{\sqrt{Q}}v$ when $X$ is an infinite dimensional
Ornstein Uhlenbeck process and $A$ and $Q$ commute.
We apply this estimate to prove that there exists a mild solution of
the semilinear Kolmogorov equation (\ref{Kolmo}) with the nonlinear term which is quadratic 
with respect to the $\sqrt{Q}$-derivative. Also in the finite dimensional case, in the setting of
\cite{BaoDeHu} and of \cite{Ri}, solution of the related Kolmogorov equation can be achieved with our techniques.

Assume that $A$ and $\sqrt{Q}$ commute and that $v$
is the unique mild solution of equation (\ref{Kolmo}).
As already noticed, if $\phi$ and $\psi$ are G\^ateaux
differentiable, then $v$ is given by $v(t,x)=Y_t^{t,x}$ where
$(Y_\tau^{t,x}, Z_\tau^{t,x})_{\tau\in[t,T]}$ solve the BSDE in the
forward-backward system (\ref{fbsde}) with $l=0$, that we rewrite here:
\begin{equation}\label{fbsde_l_nulla}
    \left\{\begin{array}{l}\dis dX_\tau =
AX_\tau d\tau+ \sqrt{Q} dW_\tau,\quad \tau\in
[t,T]\subset [0,T],
\\\dis
X_t=x,
\\\dis
 dY_\tau=-\psi(Z_\tau)\;d\tau +Z_\tau\;dW_\tau,
  \\\dis
  Y_T=\phi(X_T^{t,x}).
\end{array}\right.
\end{equation}
Moreover, by \cite{Kob}, see also \cite{BriFu}, when $\psi$
is quadratic in $Z$ it turns out that $(\Phi(\tau)=\int_t^\tau Z_sdW_s)_{\tau\in[t,T]}$
is a BMO martingale and 
\begin{equation*}
 \Vert \Phi\Vert_{\textit{BMO}}=\sup_{\sigma\in[t,T]}\E\left[ 
 \int_\sigma^T Z^2_sds \vert \calf_\sigma \right]^{1/2}<+\infty,
\end{equation*}
where the supremum is taken over all stopping times $\sigma\in[t,T]$ a.s.
Moreover, as a consequence, the stochastic exponential martingale
\begin{equation*}
\cale ( \Phi)_\tau=\cale_\tau=\exp\left(\int_t^\tau Z_sdW_s-\frac{1}{2}
 \int_t^\tau Z^2_sds  \right),
\end{equation*}
is uniformly integrable.
We are ready to prove an estimate on $Z$,
independent on $\nabla\phi$.
\begin{theorem}\label{teostimabismut}
Let $(Y,Z)$ be the solution of the BSDE in (\ref{fbsde_l_nulla}). Let $A$ and $\sqrt{Q}$
satisfy hypothesis \ref{ip su AB}, and assume that $A$ and $\sqrt{Q}$ commute. Let $\phi$ and $\psi$
satisfy hypotheses \ref{ip su psi} and \ref{ipfidiffle} with $p=2$. Then the folllowing estimate holds true:
\begin{equation}\label{stimabismut}
\vert Z_t^{t,x}\vert\leq C (T-t)^{-1/2},
\end{equation}
where $C$ depends on $t,\;T,\;A,\;\Vert\phi\Vert_\infty$ and not on $\nabla \phi$.
\end{theorem}
\dim
Let us take the $\sqrt{Q}$-derivative in the BSDE in (\ref{fbsde_l_nulla}) in the direction $h\in H$.
Let us denote $F_\tau^{t,x}=\nabla^{\sqrt{Q}}Y_\tau^{t,x}h$ and $V_\tau^{t,x}=\nabla^{\sqrt{Q}}Z_\tau^{t,x}h$.
$(F,V)$ solve the following BSDE
\begin{equation}\label{bsde_differ}
\left\{\begin{array}{l}\dis
 dF_\tau^{t,x}=-\nabla\psi(Z_\tau^{t,x})V_\tau^{t,x}\;d\tau +V_\tau^{t,x}\;dW_\tau,
  \\\dis
  F^{t,x}_T=\nabla\phi(X_T^{t,x})e^{(T-t)A}\sqrt{Q}h,
\end{array}\right.
\end{equation}
Let us denote by $\Q$ the equivalent probability measure such that
$$dW^\Q_\tau:=-\int _0^\tau \psi(Z_s^{t,x}) ds +W_\tau$$
is a Wiener process.
Notice that by our assumptions $\nabla\psi$ has linear growth
with respect to $Z$, so $\left( \int_t^\tau \nabla \psi(Z_s^{t,x})dW_s\right)_{\tau\in[t,T]}$
is a BMO martingale and
\begin{equation*}
\cale_\tau=\exp\left(\int_t^\tau \psi(Z_s^{t,x})W_s-\frac{1}{2}
 \int_t^\tau \vert\psi(Z_s^{t,x})\vert^2 ds  \right),
\end{equation*}
is a uniformly integrable martingale. Notice that
$\dfrac{d\P}{d\Q}=\cale_T$.

\noindent In $(\Omega,\calf, \Q)$, $(F^{t,x},V^{t,x})$ solve
the following BSDE:
\begin{equation}\label{bsde_differQ}
\left\{\begin{array}{l}\dis
 dF_\tau^{t,x}=V_\tau^{t,x}\;dW^\Q_\tau,
  \\\dis
  F^{t,x}_T=\nabla\phi(X_T^{t,x})e^{(T-t)A}h,
\end{array}\right.
\end{equation}
It turns out that $(F^{t,x})^2$ is a $\Q$-submartingale.

\noindent In $(\Omega,\calf, \Q)$, $(Y^{t,x},Z^{t,x})$ solve
the following BSDE:
\begin{equation*}
\left\{\begin{array}{l}\dis
 dY_\tau^{t,x}=-\psi(Z_\tau^{t,x})d\tau+\nabla\psi(Z_\tau^{t,x})Z_\tau^{t,x}\;d\tau+Z_\tau^{t,x}\;dW^\Q_\tau,
  \\\dis
  Y^{t,x}_T=\phi(X_T^{t,x}),
\end{array}\right.
\end{equation*}
and by our assumptions the generator
$-\psi(Z_\tau^{t,x})+\nabla\psi(Z_\tau^{t,x})Z_\tau^{t,x}$
has quadratic growth with respect to $Z$, so again by \cite{Kob}
$(\int_t^\tau Z_s^{t,x}ds,\tau\in[t,T])$ is a BMO $\Q$-martingale, with BMO norm depending only on $T$, $t$,
$A$ and $\Vert\phi\Vert_\infty$.

\noindent Moreover let us denote again $Y_t^{t,x}=v(t,x)$. Since $A$ and $Q$ commute,
\begin{align*}
 &F_\tau^{t,x}:=<\nabla v(\tau,X_\tau^{t,x}),\sqrt{Q}h>=
<\nabla_x v(\tau,X_\tau^{t,x}),e^{(\tau-t)A}\sqrt{Q}h> \\ \nonumber
&=<\nabla_x v(\tau,X_\tau^{t,x}),\sqrt{Q}e^{(\tau-t)A}h>=
<Z_\tau^{t,x},e^{(\tau-t)A}h>. \\ \nonumber
\end{align*}
With these facts we can prove the desired estimate, indeed, since
$(F^{t,x})^2$ is a $\Q$-submartingale
$$
\E^\Q\left[ \int_t^T \vert F_s^{t,x}\vert^2ds \vert \calf_t \right]
\geq \vert F^{t,x}_t\vert^2 (T-t)=\vert Z^{t,x}_t h\vert^2 (T-t).
$$
Moreover
$$
\E^\Q\left[ \int_t^T \vert F_s^{t,x}\vert^2ds \vert \calf_t \right]
=\E^\Q\left[ \int_t^T \vert <Z_s^{t,x},e^{(s-t)A}h>\vert^2ds \vert \calf_t \right]
\leq \Vert \Phi\Vert_{\textit{BMO},\Q},
$$
where $C_{t,T}$ depends only on $T$, $t$, $A$.
So
\begin{equation*}
\vert Z_t^{t,x}\vert\leq C (T-t)^{-1/2},
\end{equation*}
where $C$ depends on $t,\;T,\;A,\;\Vert\phi\Vert_\infty$ and not on $\nabla \phi$.
\qed

\begin{remark}\label{remarkbismut}
Notice that under the assumptions of theorem
\ref{teostimabismut} $Z_t^{t,x}=\nabla^{\sqrt{Q}}v(t,x)$, where
$v(t,x)=Y_t^{t,x}$ is the unique mild solution of equation
\ref{Kolmo}. So estimate (\ref{stimabismut})
gives a bound of $\nabla^{\sqrt{Q}}v$ in $C_{\alpha}^{s}\left(  \left[
0,T\right]  \times H,H^{\ast}\right)  $, with $\alpha=1/2$. 
\end{remark}
 Next we want to apply the result of theorem \ref{teostimabismut}
to find a mild solution to equation \ref{Kolmo}.
\begin{theorem}\label{esistenzaKolmocont}
 Let $A$ and $\sqrt{Q}$
satisfy hypothesis \ref{ip su AB} and assume that $A$ and $\sqrt{Q}$ commute. Let $\phi$ and $\psi$
satisfy hypotheses \ref{ip su psi}. Then equation (\ref{Kolmo}) admits
a unique mild solution $v$ according to definition \ref{defsolmildkolmo}.
\end{theorem}
\dim
Let $\phi\in UC_{b}\left(  H\right)  $. We can define, see e.g. \cite{DP3} and
\cite{LL}, the inf-sup convolutions $\phi_n$ of $\phi$ by setting, for $n\geq1$,
\begin{equation}
\phi_n\left(  x\right)  =\inf_{y\in H}\left\{  \phi\left(  y\right)
+2n\left\vert x-y\right\vert _{H}^{2}\right\}  , \label{inf conv}%
\end{equation}
It is well known that the inf-sup convolution $\phi_n$ of $\phi$
provides an approximation of $\phi$ in the norm of the uniform convergence, preseving the supremum norm,
and for every $n,\; \phi_n$ is lipschitz continuous and Frechet differentiable,
with derivative blowing up like $n$ as $n\rightarrow +\infty$.
So for every $n$ equation
\begin{equation}
\left\{
\begin{array}
[c]{l}%
\frac{\partial u}{\partial t}(t,x)=-\mathcal{L}u\left(  t,x\right)
+\psi\left(  \nabla^{\sqrt{Q}}u\left(  t,x\right)  \right)+l(x)  ,\text{ \ \ \ \ }t\in\left[  0,T\right]
,\text{ }x\in H\\
u(T,x)=\phi_n\left(  x\right)  ,
\end{array}
\right.  \label{Kolmo_n}%
\end{equation}
admits a unique mild solution $v_n$ according to definition
\ref{defsolmildkolmo}.
Notice that $v_n(t,x)=Y^{n,t,x}$, where we denote 
by $(X^{t,x}, Y^{n,t,x},Z^{n,t,x})$ the unique solution
of a forward backward system like (\ref{fbsde}) with final condition $\phi$
replaced by $\phi_n$.
It is immediate to see from the backward equation that
$$
\vert Y_t^{n,t,x}-Y_t^{k,t,x}\vert \leq C\Vert\phi_n-\phi_k\Vert_{\infty}
$$
Indeed, let us consider the BSDE solved by $Y^{n,t,x}-Y^{k,t,x}$
\begin{equation}\label{BSDEnk}
\left\{\begin{array}{l}\dis
 d(Y_\tau^{n,t,x}-Y_\tau^{k,t,x})=-\left(\psi(Z_\tau^{n,t,x})-\psi(Z_\tau^{k,t,x})\right)\;d\tau 
+\left(Z_\tau^{k,t,x}-Z_\tau^{k,t,x}\right)\;dW_\tau,
  \\\dis
  Y^{n,t,x}_T=\phi^n(X_T^{t,x})-\phi^k(X_T^{t,x}).
\end{array}\right.
\end{equation}
Let us denote by $\Q^{n,k}$ the equivalent probability measure such that
$$dW^{\Q^{n,k}}_\tau:=-\int _0^\tau 
\frac{\psi(Z_s^{n,t,x})-\psi(Z_s^{k,t,x})}{\vert Z_s^{n,t,x})-Z_s^{k,t,x}\vert}
\chi_{\left\lbrace Z_s^{n,t,x}-Z_s^{k,t,x}\neq 0\right\rbrace } ds +W_\tau$$
is a Wiener process.
Writing equation \ref{BSDEnk} in
$(\Omega,\calf,\Q^{n,k})$
we get the desired estimate.
\noindent Notice that by our assumptions 
$$
\frac{\psi(Z_s^{n,t,x})-\psi(Z_s^{k,t,x})}{\vert Z_s^{n,t,x}-Z_s^{k,t,x}\vert}
\chi_{\left\lbrace Z_s^{n,t,x}-Z_s^{k,t,x}\neq 0\right\rbrace}\leq C(1+\vert Z_s^{n,t,x}\vert+\vert Z_s^{k,t,x}\vert)
$$
so,
$$
\left(\int_t^\tau \frac{\psi(Z_s^{n,t,x})-\psi(Z_s^{k,t,x})}{\vert Z_s^{n,t,x}-Z_s^{k,t,x}\vert}
\chi_{\left\lbrace Z_s^{n,t,x}-Z_s^{k,t,x}\neq 0\right\rbrace}dW_s\right)_{\tau\in[t,T]}
$$
is a BMO martingale and
 \begin{align*}
 \cale^{n,k}_\tau &=\exp\left(\int_t^\tau \frac{\psi(Z_s^{n,t,x})-\psi(Z_s^{k,t,x})}{\vert Z_s^{n,t,x}-Z_s^{k,t,x}\vert}
 \chi_{\left\lbrace Z_s^{n,t,x}-Z_s^{k,t,x}\neq 0\right\rbrace}dW_s\right. \\ \nonumber
 & \left.-\frac{1}{2} \int_t^\tau \vert\frac{\psi(Z_s^{n,t,x})}-\psi(Z_s^{k,t,x}){\vert Z_s^{n,t,x}-Z_s^{k,t,x}\vert}
 \chi_{\left\lbrace Z_s^{n,t,x}-Z_s^{k,t,x}\neq 0\right\rbrace}\vert^2 ds  \right), \\ \nonumber
 \end{align*}
% \begin{equation*}
% \cale^{n,k}_\tau =\exp\left(\int_t^\tau \frac{\psi(Z_s^{n,t,x})-\psi(Z_s^{k,t,x})}{\vert Z_s^{n,t,x}-Z_s^{k,t,x}\vert}
% \chi_{\left\lbrace Z_s^{n,t,x}-Z_s^{k,t,x}\neq 0\right\rbrace}dW_s
% -\frac{1}{2} \int_t^\tau \vert\frac{\psi(Z_s^{n,t,x})-\psi(Z_s^{k,t,x})}{\vert Z_s^{n,t,x}-Z_s^{k,t,x}\vert}
% \chi_{\left\lbrace Z_s^{n,t,x}-Z_s^{k,t,x}\neq 0\right\rbrace}\vert^2 ds  \right),
% \end{equation*}
is a uniformly integrable martingale, with
$\dfrac{d\P}{d\Q^{n,k}}=\cale^{n,k}_T$.
By previous arguments we know that
$$
\vert Y_t^{n,t,x}-Y_t^{k,t,x}\vert \leq C\Vert\phi\Vert_{\infty}
$$
uniformly in $n$ and $k$.
Next we have to prove that
$$
\vert Z_t^{n,t,x}-Z_t^{k,t,x}\vert \leq C\Vert\phi_n-\phi_k\Vert_{\infty}
$$
We differentiate equation (\ref{BSDEnk}), rewritten in $(\Omega,\calf,\Q^{n,k})$,
in the direction $\sqrt{Q}h$, $h\in H$. We get
\begin{equation}\label{BSDEnk_differ}
\left\{\begin{array}{l}\dis
 d(F_\tau^{n,t,x}-F_\tau^{k,t,x})=(V_\tau^{n,t,x}-V_\tau^{k,t,x})\;dW_\tau,
  \\\dis
  F^{n,t,x}_T-F^{k,t,x}_T=(\nabla\phi^n(X_T^{t,x})-\nabla\phi^n(X_T^{t,x}))e^{(T-t)A}\sqrt{Q}h,
\end{array}\right.
\end{equation}
So, in $(\Omega,\calf,\Q^{n,k})$, $\left\lbrace F_\tau^{n,t,x}-F_\tau^{k,t,x}, \tau\in[t,T]\right\rbrace $
is a martingale:
$$
\E^{\Q^{n,k}}\left[ \int_t^T \vert F^{n,t,x}_s-F^{k,t,x}_s\vert^2
ds \vert \calf_t \right]
\geq \vert F^{n,t,x}_t-F^{k,t,x}_t\vert^2(T-t)=
\vert (Z^{n,t,x}_t-Z^{k,t,x}_t) h\vert^2  (T-t).
$$
Moreover
\begin{align*}
\E^{\Q^{n,k}}\left[ \int_t^T \vert F^{n,t,x}_s-F^{k,t,x}_s\vert^2 ds \vert \calf_t \right]&
=\E^{\Q^{n,k}}\left[ \int_t^T \vert <Z_s^{n,t,x}-Z_s^{k,t,x},e^{(s-t)A}h>\vert^2 ds \vert \calf_t \right] \\ \nonumber
&\leq C_{t,T}\E^{\Q^{n,k}}\left[ \int_t^T \vert Z_s^{n,t,x}-Z_s^{k,t,x}\vert^2 ds \vert \calf_t \right], \\ \nonumber
\end{align*}
where $C_{t,T}$ is a bounded constant depending on $t,\;T,\;A$.
Moreover in $(\Omega,\calf,\Q^{n,k})$
by equation (\ref{BSDEnk}) we immediately get
$$
\E^{\Q^{n,k}}\left[ \int_t^T \vert Z_s^{n,t,x}-Z_s^{k,t,x}\vert^2 ds \vert \calf_t \right]
\leq \Vert\phi_n-\phi_k\Vert_{\infty}
$$
and so 
\begin{equation}\label{stimaZnk}
\vert Z_t^{n,t,x}-Z_t^{k,t,x}\vert\leq C_{t,T}  \Vert\phi_n-\phi_k\Vert_{\infty}(T-t)^{-1/2},
\end{equation}
where $C_{t,T}$ is a bounded constant depending on $t,\;T,\;A$ and not on $\nabla \phi$.

\noindent So we get that by setting $v^n(t,x)=Y_t^n(t,x)$, the solution of the 
Kolmogorov equation (\ref{Kolmo_n}) $v^n(t,x)$ converges in $C([0,T]\times H)$ to $v(t,x)$, equal to
$Y_t^{t,x}$. Moreover for every $n$, $Z_t^{n,t,x}=\nabla^{\sqrt{Q}}v^n(t,x)$, and by
(\ref{stimaZnk}) $(\nabla^{\sqrt{Q}}v^n(t,x))_n$ is a Cauchy sequence in $C^s_{1/2}([0,T]\times H))$.
So $\nabla^{\sqrt{Q}}v^n(t,x)$ converges in $C^s_{1/2}([0,T]\times H))$
to an element that we denote by $F(t,x)$.
For every $n\geq 1$, 
\[
\dfrac{v^n(t,x+s\sqrt{Q}h)-v^n(t,x)}{s}\int_{0}^{1}\nabla
^{\sqrt{Q}}v^n(t,x+r\sqrt{Q}h)h \,dr.
\]
As $n\rightarrow +\infty$ we get
\[
\dfrac{v(t,x+s\sqrt{Q}h)-v(t,x)}{s}\int_{0}^{1}F(t,x+r\sqrt{Q}h)h \,dr.
\]
which gives $F(t,x)h=\nabla^{\sqrt{Q}}v(t,x)h$.
It remains to see that for every $\tau \in [0,T]$,
$\nabla^{\sqrt{Q}}v(\tau,X_\tau^{t,x})h=Z_\tau^{t,x}$,
where $Z^{t,x}$ is the limit of $Z^{n,t,x}$ in
$L^2(\Omega\times[0,T])$. It turns out that by previous calculations
$\nabla^{\sqrt{Q}}v^n(\tau,X_\tau^{t,x})\rightarrow
\nabla^{\sqrt{Q}}v(\tau,X_\tau^{t,x})$ in $C^s_{1/2}([0,T]\times H))$. So
$\nabla^{\sqrt{Q}}v(\tau,X_\tau^{t,x})=Z_\tau^{t,x}$ $\P$ a.s. for a.a.
$\tau\in[t,T]$.
Since $(Y,Z)$ solve the BSDE in (\ref{fbsde_l_nulla}), with $Y_t^{t,x}=v(t,x)$,
by previous arguments we get $Z_t^{t,x}=\nabla^{\sqrt{Q}}v(t,x)$.
By classical arguments we deduce that $v$ solves equation (\ref{Kolmo}).
Moreover the solution is unique since the solution of the corresponding BSDE
is unique.
\qed

\subsection{A quadratic optimal control problem}
\label{appl-contr-2}
We apply the result of the previous section to a control problem where the current
cost has quadratic growth with respect to the control $u$ and the final cost is only continuous.
The fundamental relation and the existence of a solution of the closed loop 
equation cannot be achieved as in theorem \ref{th-rel-font} and \ref{teo su controllo debole}
respectively, since this time $\Z^{t,x}_t=\nabla^{\sqrt{Q}}v(t,x)$ is not bounded.

Let $X^u$ the solution of equation (\ref{sdecontrolforte}), and we have to minimize
the cost functional (\ref{cost}) over all the admissible control $u$,
where by admissible control we mean here an $(\calf_t)_t$-predictable
process, taking values in a closed subset $K$ of a normed space $U$, such that
\[
 \E\int_0^T \vert u_s \vert ^2 ds < + \infty.
\]
This assumptions is natural this time since we assume here that the cost 
has quadratic growth at infinity, namely the cost must satisfy hypothesis
\ref{ip costo} with $q=2$. We define the hamiltonian function in a classical way as in
\ref{hamilton}. The hamiltonian satisfies the properties stated in lemma \ref{lemma-hamilton}, in particular
estimates (\ref{stima1psi}) and (\ref{stima2psi}) hold true with $p=2$ and the
infimum is achieved in a ball of radius $C(1+\vert z\vert)$.
\begin{theorem}\label{th-rel-font-quadr}
 Assume that $A$ and $\sqrt{Q}$ satisfy hypothesis \ref{ip su AB} and commute.
Let $g$ and $l$ satisfy point 2 and 3 of hypothesis \ref{ip costo}
and let hypothesis \ref{ip aggiuntive} hold true; assume that the hamiltonian function $\psi$
satisfies G\^ateaux differentiability assumptions stated in hypothesis \ref{ip su psi}.
 For
every $t\in [0,T]$, $x\in H$ and for
 all admissible control $u$ we have $J(t,x,u(\cdot))
 \geq v(t,x)$,
  and the
 equality holds if and only if
$$
u_s\in \Gamma\left( \nabla^{\sqrt{Q}}
v(s ,X^{u,t,x}_s)
\right)
  $$
\end{theorem}
{\bf Proof.} The proof follows from proposition 4.1 in \cite{fuhute},
and by our assumption here we have also the identification
$Z_\tau^{t,x}=\nabla^{\sqrt{Q}}v(\tau, X_\tau^{t,x}).$

With the assumptions of Theorem
\ref{th-rel-font-quadr}, we can define the so called
optimal feedback law as we have done in \ref{leggecontrolloottima}.
Since as we have already noticed in section \ref{applic contr 1} existence of a solution of the closed loop
equation is not obvious, we formulate the optimal control problem
in the weak sense, following \cite{FlSo}, see also section \ref{applic contr 1}.

\begin{theorem}\label{teo su controllo debole quadr}
Assumehat $A$ and $\sqrt{Q}$ satisfy hypothesis \ref{ip su AB} and commute.
Let $g$ and $l$ satisfy point 2 and 3 of hypothesis \ref{ip costo}
and let hypothesis \ref{ip aggiuntive} hold true; assume that the hamiltonian function $\psi$
satisfies G\^ateaux differentiability assumptions
stated in hypothesis \ref{ip su psi}.
 For
every $t\in [0,T]$, $x\in H$ and for
 all admissible control systems we have $J(t,x,u(\cdot))
 \geq v(t,x)$,
  and the
 equality holds if and only if
$$
u_s\in \Gamma\left( \nabla^{\sqrt{Q}}
v(s ,X^{u}_s)
\right)
  $$

Moreover
assume that the set-valued map $\Gamma$ and let $\gamma$
be its measurable selection. 
\begin{equation*}
u_{\tau}=\gamma(\nabla^{\sqrt{Q}} v(\tau,X^u_\tau ))
,\text{ \ }\mathbb{P}\text{-a.s. for a.a. }\tau
\in\left[ t,T\right]
\end{equation*}
is optimal.

Finally, the closed loop equation \ref{closed loop eq}
admits a weak solution
$(\Omega,\mathcal{F},
\left(\mathcal{F}_{t}\right) _{t\geq 0}, \mathbb{P}, W,X)$
which is unique in law and setting
$$
u_{\tau}=\gamma\left(\nabla^{\sqrt{Q}} v(\tau,X_\tau )\right),
$$
we obtain an optimal admissible control system $\left(
W,u,X\right) $.
\end{theorem}

\noindent {\bf Proof.} The proof follows from the fundamental relation
stated in theorem \ref{th-rel-font-quadr}.
For the solvability of the closed loop equation we refer to proposition
5.2 in \cite{fuhute}.
\qed

\section{Optimal control problems for the heat equation}
\label{sez_contr_heat}

In this section we present some control problem related to a stochastic heat equation.
As in section \ref{sezionesde}, when introducing equation (\ref{heat equation}),
here $\calo$ is a bounded domain in $\R$,
$H=L^2(\calo)$ and $\left\lbrace e_k\right\rbrace_{k\in\N}$
is the complete orthonormal basis which diagonalizes $\Delta$,
endowed with Dirirchlet boundary conditions in $\calo$.
$Q:H\rightarrow H$ satisfies hypothesis \ref{ip-cov}, 
in particular $Q e_k =\lambda_k e_k,\,\lambda_k\geq0, \;k\in\N$.
We consider the following controlled heat equation
\begin{equation}\label{heat equation contr}
 \left\{
  \begin{array}{l}
  \dis
\frac{ \partial y}{\partial s}(s,\xi)= \Delta y(s,\xi)+
+\sum_{k\in\N}\sqrt{\lambda_k} \left( \int_{\calo} u_s(\eta)e_k(\eta)d\eta\right)
e_k(\xi)+\frac{ \partial W^Q
}{\partial s}(s,\xi), \qquad s\in [t,T],\;
\xi\in \calo,
\\\dis
y(t,\xi)=x(\xi),
\\\dis
 y(s,\xi)=0, \quad \xi\in\partial \calo.
\end{array}
\right.
\end{equation}
where $u_s\in L^2(\calo)$ represents the control.
In the following we denote by 
$\mathcal{A}_{d}$ the set of
admissible controls, that is the real valued predictable processes such that 
\[
 \E \int_0^T \left( \int_{\calo}\vert u_t(\xi) \vert^2 d\xi\right)  ^{q/2} dt <+\infty.
\]
and such that $u_t \in K$, where $K$ is a closed subset of $H$, not
necessarily coinciding with $H$.
Here as in equation (\ref{heat equation}),  $W^Q(s,\xi)$ is a Gaussian mean zero random field, such that
the operator $Q$ characterizes the correlation in the space variables.
Our aim is to minimize over all admissible controls the cost functional
\begin{equation}
J\left(  t,x(\xi),u\right)  =\mathbb{E}\int_{t}^{T}\int_{\calo}[
\bar l\left(X^{u}_s(\xi)\right)+\vert u_s(\xi)\vert^q ]d\xi ds
+\mathbb{E}\int_{\calo} \bar\phi\left(X^{u}_T(\xi)\right)\,d\xi. 
\label{heat cost}%
\end{equation}
for real functions $\bar\phi$ and $\bar l$, and for $q\leq 2$.

We make the following assumptions on the cost $J$.

\begin{hypothesis}
\label{ip costo heat}
The function $\mathbb{\bar\phi}:\R\rightarrow\mathbb{R}$ is lipschitz
continuous and bounded; $\bar l: \R\rightarrow\mathbb{R}$ is bounded and continuous.
% $g:\R\rightarrow\R$ is continuous and for every $u\in\R$, for $1<q\leq 2$
% \begin{equation*}
% 0\leq g(u)\leq c(1+ \vert u \vert ^q) 
% \end{equation*}
% and there exist $R>0$, $C>0$ such that
% \begin{equation*}
% g(u)\geq C \vert u \vert ^q \qquad \text{for every }u \in K : \vert u\vert \geq R.
% \end{equation*}
\end{hypothesis}
Let us define, for $\xi\in H$
\begin{equation}\label{heat notazioni}
 \phi (x)=\int_{\calo} \bar\phi\left(x(\xi)\right)\,d\xi,\quad
l(x)=\int_{\calo}
\bar l\left(x(\xi)\right)\,d\xi.
\end{equation}
It turns out that if $\bar l$ and $\bar\phi$ satisfy
hypothesis \ref{ip costo heat}, then $\phi$ and $l$
defined in (\ref{heat notazioni}) satisfy hypothesis
\ref{ip costo}.
Moreover by defining $g(u)=\int_\calo\vert u_s(\xi)\vert^q d\xi=\vert u\vert_{L^2(\calo)}$,
then the hamiltonian function turns out to be
$\psi(z)= (\frac{1}{q})^{\frac{1}{q-1}}\frac{1-q}{q}\vert z\vert ^p$.

\noindent Moreover equation (\ref{heat equation contr}) can be written in an abstract way in $H$
as
\begin{equation}\label{heat eq abstr contr}
 \left\{
\begin{array}
[c]{l}%
dX_\tau  =AX_\tau d\tau+\sqrt{Q}u_\tau+\sqrt{Q}dW_\tau
,\text{ \ \ \ }\tau\in\left[  t,T\right] \\
X_t =x,
\end{array}
\right.
\end{equation}
where $A$ is the Laplace operator with Dirirchlet boundary conditions,
$W$ is a cylindrical Wiener process in $H$ and $Q$ is its covariance operator.
The control problem in its abstract formulation is to minimize
over all admissible controls the cost functional
\begin{equation}
J\left(  t,x,u\right)  =\mathbb{E}\int_{t}^{T} [l\left(X^{u}_s\right)
+\vert u_s\vert^q  ]ds
+\mathbb{E}\phi\left(X^{u}_T\right).
\label{heat cost abstr}%
\end{equation}
By applying results in section \ref{applic contr 1}, we get the following

\begin{theorem}\label{heat-th-contr}
Let $X^u$ be the solution of equation (\ref{heat equation contr}), let the cost
be defined as in (\ref{heat cost}) and let
\ref{ip costo heat} hold true. Moreover assume
that the hamiltonian function $\psi$
satisfies G\^ateaux differentiability assumptions
stated in hypothesis \ref{ip su psi}.
 For every $t\in [0,T]$, $x\in L^2(\calo)$ and for
 all admissible control $u$ we have $J(t,x,u(\cdot))
 \geq v(t,x)$,
  and the
 equality holds if and only if
$$
u_s\in \Gamma\left( \nabla^{\sqrt{Q}}
v(s ,X^{u,t,x}_s)
\right)
  $$
Moreover
assume that the set-valued map $\Gamma$ is nonempty and let $\gamma$
% be its measurable selection. 
% \begin{equation*}
% u_{\tau}=\gamma(\nabla^{\sqrt{Q}} v(\tau,X^u_\tau ))
% ,\text{ \ }\mathbb{P}\text{-a.s. for a.a. }\tau
% \in\left[ t,T\right]
% \end{equation*}
% is optimal.

The closed loop equation admits a weak solution
$(\Omega,\mathcal{F},
\left(\mathcal{F}_{t}\right) _{t\geq 0}, \mathbb{P}, W,X)$
which is unique in law and setting
$$
u_{\tau}=\gamma\left(\nabla^{\sqrt{Q}} v(\tau,X_\tau )\right),
$$
we obtain an optimal admissible control system $\left(
W,u,X\right) $.
\end{theorem}

\noindent {\bf Proof.} The proof follows from the abstract formulation of the problem, and 
by applying theorems \ref{th-rel-font} and \ref{teo su controllo debole}.
\qed

Next we turn to an optimal control problem related to the
controlled equation
(\ref{heat eq abstr contr}) with quadratic cost $g$ and consequently
quadratic hamiltonian function, and with final cost continuous.
In this case, in order to perform the synthesis of the optimal control,
we apply the results of section \ref{appl-contr-2}.
Namely we consider equation (\ref{heat equation contr}). We have to minimize
the cost functional
\begin{equation}
J\left(  t,x(\xi,u\right)  =\mathbb{E}\int_{t}^{T}\int_{\calo}[
\bar l\left(X^{u}_s(\xi)\right)+\bar g( u_s(\xi)) ]d\xi ds
+\mathbb{E}\int_{\calo} \bar\phi\left(X^{u}_T(\xi)\right)\,d\xi. 
\label{heat cost quadr}%
\end{equation}
over all admissible controls, that is real valued predictable processes such that 
\[
 \E \int_0^T \left( \int_{\calo}\vert u_t(\xi) \vert^2 d\xi\right)   dt <+\infty.
\]
$\bar\phi$, $\bar g$ and $\bar l$ are real functions satisfying the following:
\begin{hypothesis}
\label{ip costo heat quadr}
The function $\mathbb{\bar\phi}:\R\rightarrow\mathbb{R}$ is 
continuous and bounded;  $\bar l: \R\rightarrow\mathbb{R}$ is bounded and continuous;
$g:\R\rightarrow\R$ is continuous and for every $u\in\R$
\begin{equation*}
0\leq g(u)\leq c(1+ \vert u \vert ^2) 
\end{equation*}
and there exist $R>0$, $C>0$ such that
\begin{equation*}
g(u)\geq C \vert u \vert ^2 \qquad \text{for every }u \in K : \vert u\vert \geq R.
\end{equation*}
\end{hypothesis}
Equation (\ref{heat equation contr}) admits the abstract formulation given by
(\ref{heat eq abstr contr}) and the cost functional can be formulated in an abstract way
as
\begin{equation*}
J\left(  t,x,u\right)  =\mathbb{E}\int_{t}^{T} l\left(X^{u}_s\right)
+g(u_s) ds
+\mathbb{E}\phi\left(X^{u}_T\right).
\end{equation*}
with notation \ref{heat notazioni} and by setting moreover
\begin{equation*}
 g (u)=\int_{\calo} \bar g\left(u(\xi)\right)\,d\xi,\quad.
\end{equation*}
It turns out that if $\bar\phi$, $\bar l$ and $\bar g$
satisfy hypothesis \ref{ip costo heat quadr}, then $\phi$,
$l$ and $g$ satisfy hypothesis \ref{ip costo} with $q=2$.

\begin{theorem}\label{heat-th-contr-quadr}
Let $X^u$ be the solution of equation (\ref{heat equation contr}), let the cost
be defined as in \ref{heat cost quadr} and let
\ref{ip costo heat quadr} hold true. Moreover assume
that the hamiltonian function $\psi$
satisfies G\^ateaux differentiability assumptions
stated in hypothesis \ref{ip su psi}.
 For every $t\in [0,T]$, $x\in L^2(\calo)$ and for
 all admissible control $u$ we have $J(t,x,u(\cdot))
 \geq v(t,x)$,
  and the
 equality holds if and only if
$$
u_s\in \Gamma\left( \nabla^{\sqrt{Q}}
v(s ,X^{u,t,x}_s)
\right)
  $$
Moreover
assume that the set-valued map $\Gamma$ is nonempty and let $\gamma$
be its measurable selection. The closed loop equation
admits a weak solution
$(\Omega,\mathcal{F},
\left(\mathcal{F}_{t}\right) _{t\geq 0}, \mathbb{P}, W,X)$
which is unique in law and setting
$$
u_{\tau}=\gamma\left(\nabla^{\sqrt{Q}} v(\tau,X_\tau )\right),
$$
we obtain an optimal admissible control system $\left(
W,u,X\right) $.
\end{theorem}

\noindent {\bf Proof.} The proof follows from the abstract formulation of the problem, and 
by applying theorem \ref{th-rel-font-quadr}.
\qed


\begin{thebibliography}{11}




 \bibitem{AuFr} J.P. Aubin, H. Frankowska, Set valued analysis,
 Birkh\"{a}user, Boston, 1990.
 

\bibitem{BaoDeHu}  X. Bao, F. Delbaen and Y. Hu,
Backward SDEs with superquadratic growth, arXiv:0902.3316,
to appear on Probability Theory and Related Fields, 
DOI: 10.1007/s00440-010-0271-1.

% \bibitem{barbu} V. Barbu, Analysis and control of nonlinear
% infinite dimensional systems, Academic Press Inc., San Diego,
% 1993.
\bibitem{BriFu} P. Briand, F. Confortola, BSDEs with stochastic Lipschitz condition and quadratic PDEs in Hilbert spaces.  Stochastic Process. Appl.  118  (2008),  no. 5, 818--838.

\bibitem{BriHu} P. Briand, Y. Hu, 
Stability of BSDEs with random terminal time 
and homogenization of semilinear elliptic PDEs.
J. Funct. Anal.  155  (1998),  no. 2, 455--494


\bibitem{Ce}S. Cerrai, A Hille-Yosida theorem for weakly continuous semigroups.
Semigroup Forum  49  (1994),  no. 3, 349--367.

\bibitem{CeGo}S. Cerrai and F. Gozzi, Strong solutions of Cauchy problems
associated to weakly continuous semigroups.  Differential Integral Equations  8  (1995),  no. 3, 465--486.

\bibitem {DP1}G. Da Prato and J. Zabczyk, Stochastic equations in
infinite dimensions, Encyclopedia of Mathematics and its Applications 44,
Cambridge University Press, 1992.


%\bibitem{daza2} G. Da Prato, J. Zabczyk,
% Ergodicity for infinite-dimensional systems.
% London Mathematical Society Lecture Note Series, 229.
% Cambridge University Press, Cambridge, 1996.
% 
% 
\bibitem {DP3}G. Da Prato and J. Zabczyk,Second order partial
 differential equations in Hilbert spaces. London Mathematical Society Note
 Series, 293, Cambridge University Press, Cambridge, 2002.

% 
% 
% \bibitem{ElMa} N. El Karoui, L. Mazliak ed.:
%  Backward Stochastic Differential Equations, Pitman Research Notes in Mathematics Series
%  364, Longman, 1997.
% 
% \bibitem{kapequ} N. El Karoui, S. Peng, M. C. Quenez, Backward
% stochastic differential equations in finance. Mathematical Finance
% 7 (1997), no. 1, 1-71.


\bibitem{FlSo} W. H. Fleming, H. M. Soner,
Controlled Markov processes and viscosity solutions. Applications
of Mathematics 25. Springer-Verlag, 1993.


\bibitem{fu} M. Fuhrman, Smoothing properties of nonlinear stochastic equations in Hilbert spaces.  NoDEA Nonlinear Differential Equations Appl.  3  (1996),  no. 4, 445--464.


\bibitem{fuhute} M. Fuhrman, Y. Hu, G. Tessitore, On a class of stochastic
optimal control problems related to BSDEs with quadratic growth. 
SIAM J. Control Optim.  45  (2006),  no. 4, 1279--1296.

\bibitem{fute} M. Fuhrman, G. Tessitore, Nonlinear Kolmogorov
equations in infinite dimensional spaces: the backward stochastic
differential equations approach and applications to optimal
control. Ann. Probab. 30 (2002), no. 3, 1397--1465.




% 
% \bibitem{fute2} M. Fuhrman, G. Tessitore, Infinite horizon
% backward stochastic differential equations and elliptic equations
% in Hilbert spaces. Ann. Probab. 32 (2004), no. 1B, 607--660.
% 
% 
% \bibitem{fute3} M. Fuhrman, G. Tessitore.
% Generalized directional gradients, backward stochastic differential equations
% and mild solutions of semilinear parabolic equations.
% Appl. Math. Optim. 51 (2005), no. 3, 279--332.


\bibitem{Go1} F. Gozzi, Regularity of solutions of second order
Hamilton-Jacobi equations in Hilbert spaces and applications to a control
problem, (1995) Comm. Partial Differential Equations 20, pp. 775-826.

\bibitem{Go2} F. Gozzi, Global regular solutions of second order
Hamilton-Jacobi equations in Hilbert spaces with locally Lipschitz
nonlinearities, (1996) J. Math. Anal. Appl. 198, pp. 399-443.

%\bibitem{ImdoRe} P. Imkeller, G. do Reis, spa 2010.

\bibitem{Kob} M. Kobylanski, Backward stochastic differential equations and partial differential equations with quadratic growth.  Ann. Probab.  28  (2000),  no. 2, 558--602.

\bibitem{LL}J. M. Lasry and P. L. Lions, A remark on regularization in
Hilbert spaces, (1986) Israel. J. Math. 55, pp. 257-266.

\bibitem{Mas} F. Masiero, Semilinear Kolmogorov equations and
applications to stochastic optimal control, Appl. Math. Optim.,
51 (2005), pp.~201--250.


 \bibitem{Mas1} F. Masiero, Regularizing properties for transition semigroups
 and semilinear parabolic equations in Banach spaces.
 Electron. J. Probab.  12  (2007), no. 13, 387--419

%\bibitem{Nu} D. Nualart,
% The Malliavin calculus and related topics.
% Probability and its Applications. Springer-Verlag, New York, 1995.


%\bibitem{NuPa} D. Nualart, E. Pardoux,
% Stochastic calculus with anticipative integrands. Probab.
% Th. Rel. Fields 78 (1988), 535-581.


\bibitem{Pa} E. Pardoux, BSDE's, weak convergence and homogeneization
of semilinear PDE's. In: Nonlinear analysis, differential
equations and control, eds. F.H. Clarke, R.J. Stern, 503-549,
Kluwer, 1999.

\bibitem{PaPe1} E. Pardoux, S. Peng,
Adapted solution of a backward stochastic differential
equation, Systems and Control Lett. 14, 1990, 55-61.

\bibitem{PaPe} E. Pardoux, S. Peng,
Backward stochastic differential equations and quasilinear
parabolic partial differential equations, in: Stochastic
partial differential equations and their applications, eds. B.L.
Rozowskii, R.B. Sowers, 200-217, Lecture Notes in Control Inf.
Sci. 176, Springer, 1992.


\bibitem{Ri} A. Richou, Numerical simulation of BSDEs with drivers of quadratic growth,
 arXiv:1001.0401.
% 
% \bibitem{ruva1} F. Russo, P. Vallois,
% Forward, backward and symmetric stochastic integration.
% Probab. Theory Related Fields 97 (1993), no. 3, 403--421.
% 
% \bibitem{ruva2} F. Russo, P. Vallois,
%  The generalized covariation process and Ito formula.
%  Stochastic Process. Appl. 59 (1995), no. 1, 81--104.
% 
% \bibitem{ruva3} F. Russo, P. Vallois,
%  Ito formula for $C\sp 1$-functions of semimartingales.
%  Probab. Theory Related Fields 104 (1996), no. 1, 27--41.
% 
% \bibitem{ruva4} F. Russo, P. Vallois,
%  Stochastic calculus with respect to continuous finite
%  quadratic variation processes.
%  Stochastics Stochastics Rep. 70 (2000), no. 1-2, 1--40.
% 



\end{thebibliography}
\end{document}